\documentclass[preprint,12pt]{elsarticle}

\usepackage[pdftex]{hyperref}
\biboptions{sort&compress,numbers,square}

\usepackage{amsfonts}
\usepackage{graphicx}
\usepackage{epsfig}
\usepackage{geometry}
\usepackage{tikz}
\usetikzlibrary{trees,positioning,shapes}
\usepackage{subcaption}
\usepackage{amsmath}
\usepackage{amsthm}

\usepackage{xcolor, colortbl}

\usepackage{algorithm}
\usepackage{algorithmicx, algpseudocode}

\newcommand{\revision}[1]{\textcolor{black}{#1}}
\newcommand{\inblue}[1]{\textcolor{blue}{#1}}
\newcommand{\inred}[1]{\textcolor{red}{#1}}

\newcommand{\Aniarev}[1]{\textcolor{black}{#1}}

\usepackage{bm}
\usepackage{mathtools}

\usepackage{rotating}

\begin{document}

\begin{frontmatter}

\title{Complexity of direct and iterative solvers on space-time formulations versus time--marching schemes for h-refined grids towards singularities}

\author{Marcin Skotniczny$^{(1)}$, Anna Paszy\'nska$^{(2)}$, Sergio Rojas$^{(3)}$, Maciej Paszy\'{n}ski$^{(1)}$}

\address{$^{(1)}$AGH University of Science and Technology,
Krak\'{o}w, Poland \\
e-mail: skotniczny.marcin@gmail.com,  maciej.paszynski@agh.edu.pl\\
$^{(2)}$ Jagiellonian University, Krak\'ow, Poland \\
e-mail: anna.paszynska@uj.edu.pl \\
$^{(3)}$ Instituto de Matem\'aticas. Pontificia Universidad Cat\'olica de Valpara\'iso, Chile \\
e-mail: sergio.rojas.h@pucv.cl}

\begin{abstract}
We study computational complexity aspects for Finite Element formulations considering hypercubic space--time full and time--marching discretization schemes for $h$--refined grids towards singularities.
We perform a relatively comprehensive study of comparing the computational time via time complexities of direct and iterative solvers.
We focus on the space-time formulation with refined computational grids and on the corresponding time slabs, namely, computational grids obtained by taking the ``cross-sections" of the refined space-time mesh.
We compare the computational complexity of the space-time formulation and the corresponding time--marching scheme.
Our consideration concerns the computational complexity of the multi-frontal solvers, the iterative solvers, as well as the static condensation.
Numerical experiments with Octave confirm our theoretical findings.
\end{abstract}
	

\end{frontmatter}

\section{Introduction}
Several approaches exist in the Finite Element (FE) community to deal with time-dependent Partial Differential Equations (PDEs). Nowadays, two common approaches are classical time--marching FE schemes (see e.g.~\cite{thomee07}) and space-time formulations, where time is treated as another spatial variable (see e.g.~\cite{langer2019}). Both methodologies enjoy several advantages, implying in general that one can run with some benefit concerning the other depending on the result of interest. For instance, if only the final state is of interest, time--marching schemes are typically preferred due to their flexibility in the election for schemes to be considered in the time discretization procedure and their optimal storage requirements as a consequence of their sequential nature. However, when considering problems where the continuous history of the evolution problem is of interest, typically, a space-time formulation will be preferred, despite the natural increment of the complexity associated with its resolution due to the dimensional increment in the discrete formulation. Moreover, space-time formulations run with an advantage in the presence of space--time singularities since they allow for the construction of robust a posteriori error estimates for local space-time mesh refinements. Time--marching schemes must remain with the sequential nature of the refinement, implying a possible increment in the storage requirements depending on the nature of the singularity. Finally, space-time formulations allow for parallel static condensation, thus possibly being faster than time--marching schemes when a sufficient amount of computational power is available. Even though there is still extensive research on the space-time formulations and the time--marching schemes, to the authors' knowledge, there are no works related to theoretical aspects of the computational complexity of both methodologies available in the literature.\\
In this work, we study complexity aspects of time--marching and space-time formulations, being the first step to further work in complexity aspects. We explore the complexity of iterative and direct solvers when considering hypercubic $h$--refined grids towards singularities. For this purpose, we propose a general method-independent strategy simulating the best-scenario possibility for each case. We proceed as follows: To simulate an $h$--adaptive space--time procedure, we start from a uniform $d$--dimensional mesh, where $d-1$ dimensions correspond to the spatial discretization and the last dimension to the time discretization, and we consider regular hypercubic geometrical refinements towards singularities. To simulate a time--marching scheme, we consider a sequence of $d-1$ dimensional meshes obtained from considering time cuts of the space--time grids, with a fixed time step of length equal to \Aniarev{the smallest element dimension in the time axis direction}. To motivate these assumptions, the reader can consider a parabolic problem, for instance, the heat equation, as the PDE of interest.
Therefore, the assumption for the space--time formulation is natural, assuming that information of the singularity is available, for instance, an a posteriori error estimator. In contrast, the assumption for the time--marching scheme will be the ideal scenario for the explicit-in-time Euler scheme ensuring that the \Aniarev{Courant–Friedrichs–Lewy condition} (CFL condition)~\cite{CFL} is satisfied, and also information on the singularity is provided.\\
To derive the estimations, we start by considering refinements toward a space-time ``edge'', resulting from a point traveling through space and time. This space-time refinement pattern corresponds to a sequence of spatial meshes refined towards ``points'' located on the space-time edge at particular time moments. Next, we consider the refinement towards a space-time ``face'', resulting from an edge traveling through space and time. Again, this space-time refinement pattern corresponds to a sequence of spatial meshes refined towards ``edges'' located on the space-time face at a particular moment.
Finally, we consider the refinement towards a space-time ``hyperface'' resulting from a face traveling through space and time. This space-time refinement pattern corresponds to a sequence of spatial meshes refined towards ``faces'' located on the space-time hyperface at a particular time instant.
We compare the computational complexity of the iterative solver executed over the $d$ dimensional space-time domain to the computational complexity of the iterative solver performed multiple times over the $d$-1 dimensional meshes during the time--marching scheme. We compare the computational complexity of the direct solver run over the $d$ dimensional space-time domain to the computational complexity of the direct solver executed multiple times over $d$-1 dimensional meshes within the time--marching scheme. \\
%We show that the time complexity of space-time formulations is usually higher than the time complexity of time marching schemes.\\
%
We summarize the findings of this work in the following two points:
\begin{itemize}
\item This paper derives a general formula for the number of unknowns for the mesh of arbitrary dimension $d$, refined towards a singularity of \Aniarev {dimension $q$}.
This paper generalizes the results discussed in \cite{ICCS2020}, \Aniarev{for the case of point singularity ($q=0$), edge singularity ($q=1$), face singularity ($q=2$), hyperface singularity ($q=3$)}, focusing on the computational complexities of space-time formulations and time--marching schemes. In our general estimates, we do not consider the polynomial order of approximation $p$, and we assume that this is a constant in our formulas.
Additionally, the singularities in our case may have an arbitrary shape, as is presented in Figure \ref{ch5:edgenonregular}. Therefore, we only fix the dimension of the space $d$ and the \Aniarev {dimension} of the singularity $q$.
\item This paper estimates the time complexity of finite element method simulations performed on adaptive $d$--dimensional space--time meshes. We also assess the time complexity of finite element method simulations performed on a sequence of adaptive $d-1$--dimensional meshes resulting from the corresponding time--marching scheme.
In particular, we estimate $N$ = the number of unknowns for the whole adaptive space--time formulation, $n$= the number of unknowns from a mesh from a sequence of adaptive meshes resulting from the time--marching scheme.
We estimate the time complexity of the iterative and direct solvers.
Thus, we have the lower bounds (expressed by the number of unknowns) and upper bounds (described by the time complexity of the sparse direct solver).
%\Aniarev{CHYBA TO ZDANIE JEST NADMIAROWE We also compare the time complexity of iterative solvers executed for either space--time formulation or time--marching scheme.}
\end{itemize}
%
%\Aniarev{NA KONCU ZMIENIC STRUCTURE OF THE WORK}
The structure of the paper is the following. We start Section 2 with some preliminaries, including a general algorithm refining the mesh towards a given singularity, followed by an overview of the sparse Gaussian elimination, usually implemented using a multi-frontal solver approach, and our notion of element partition tree that allows constructing an ordering for sparse matrix permutation in the space--time formulations. We also discuss the recursive formula for estimating the time complexity of our sparse direct solvers based on arrangements constructed from the element partition trees.
Section 3 summarizes our findings in the context of space--time formulations and time--marching schemes. We also present numerical verification using Octave codes in Section 4. The paper is concluded in the last section.
In Appendix A, we derive time complexity for a multi-dimensional grid refined toward point singularity. In Appendix B, we derive time complexity for multi-dimensional grids refined toward arbitrary singularities. In this case, we also include the dependency in the polynomial approximation orders.

\section{Preliminaries}
\subsection{Complexity literature review and assumptions}
Computational complexity, especially time complexity, is one of the most fundamental concepts of theoretical computer science. It was first defined in 1965 by Hartmanis and Stearns \cite{HS65}. \Aniarev{In this paper we estimate the computational complexity of solving systems of linear equations for multidimensional meshes with arbitrary dimensional similarities encountered in space --time formulation as well as time--marching schemes.}
The systems of equations generated by non-stationary problems can be solved by either direct \cite{DS1,DS2} or iterative solvers \cite{Saad}. The time complexity of iterative solvers, in general, can be estimated as $n_{iter} \times N$, where $n_{iter}$ is the number of iterations, and in general, it depends on the spectral properties of the matrix, and it grows with the problem size $N$.
The time complexity of direct solvers \cite{DS1,DS2} for certain classes of meshes, especially regular meshes, is well known.
In particular, for three-dimensional uniform grids, the computational complexity is of the order of ${\cal O}(N^2)$. For the two-dimensional grids, the complexity is of the order of ${\cal O}(N^{1.5})$ \cite{Liu,Singularities1}.
The sparse direct solvers rely on the smart permutation of the matrix, resulting in its banded structure and efficient sparse factorization, avoiding zero entries.
The problem of finding an optimal order of elimination of unknowns for the direct solver, in general, is indeed NP-complete \cite{NP}.
There are several heuristic algorithms analyzing the sparsity pattern of the resulting matrix
\cite{c21,c22,c23,c24,c25}.
We focus on adaptive grids, refined in space--time and space domains.
For three-dimensional grids adapted towards the point, edge, and face, the computational complexities are ${\cal O}(N), {\cal O}(N)$, and ${\cal O}(N^{1.5})$, respectively \cite{Singularities2}.
These estimates assume a prescribed order of eliminating variables \cite{threeD1}.
Similarly, for two-dimensional grids refined towards a point or edge, it is ${\cal O}(N)$ \cite{twoD1}. Again, these estimates assume a prescribed order of elimination of variables \cite{twoD2}. \Aniarev{We generalized the results of \cite{Singularities2} into multidimensional adaptive grids with arbitrary dimensional singularities.}
There are also solvers based on hierarchical matrices decomposition \cite{Hackbush}. Their time complexity generally depends on the number of unknowns in the mesh. However, the constant in front of the complexity also grows with the mesh dimension, sparsity of the matrix and rank of sub-blocks. We leave their analysis in the space--time set up for our future work.
%\Aniarev{tu maciek cos miał dodać-kartka}
%\inred{In our estimates, we assume that there is one degree of freedom per element. This assumption is justified because complexity of the local element assemble is proportional to some power of the polynomial degree $p$, that we considered fixed. Fix me babe. There is still an extensive research going on the space--time formulations and the time--marching schemes, and this paper allows to compare the computational complexity of both methods. Our future work may involve detailed comparison of modern variations of space--time formulations including time-slabs, tent pitching methods among others \cite{timeslabs,tentpitching}. Move me babe.}.

\subsection{Hierarchical basis functions and construction of an $h$-adaptive mesh}\label{ch3:algorithm-props}
\begin{figure}
    \centering
    \includegraphics[width=0.9\textwidth]{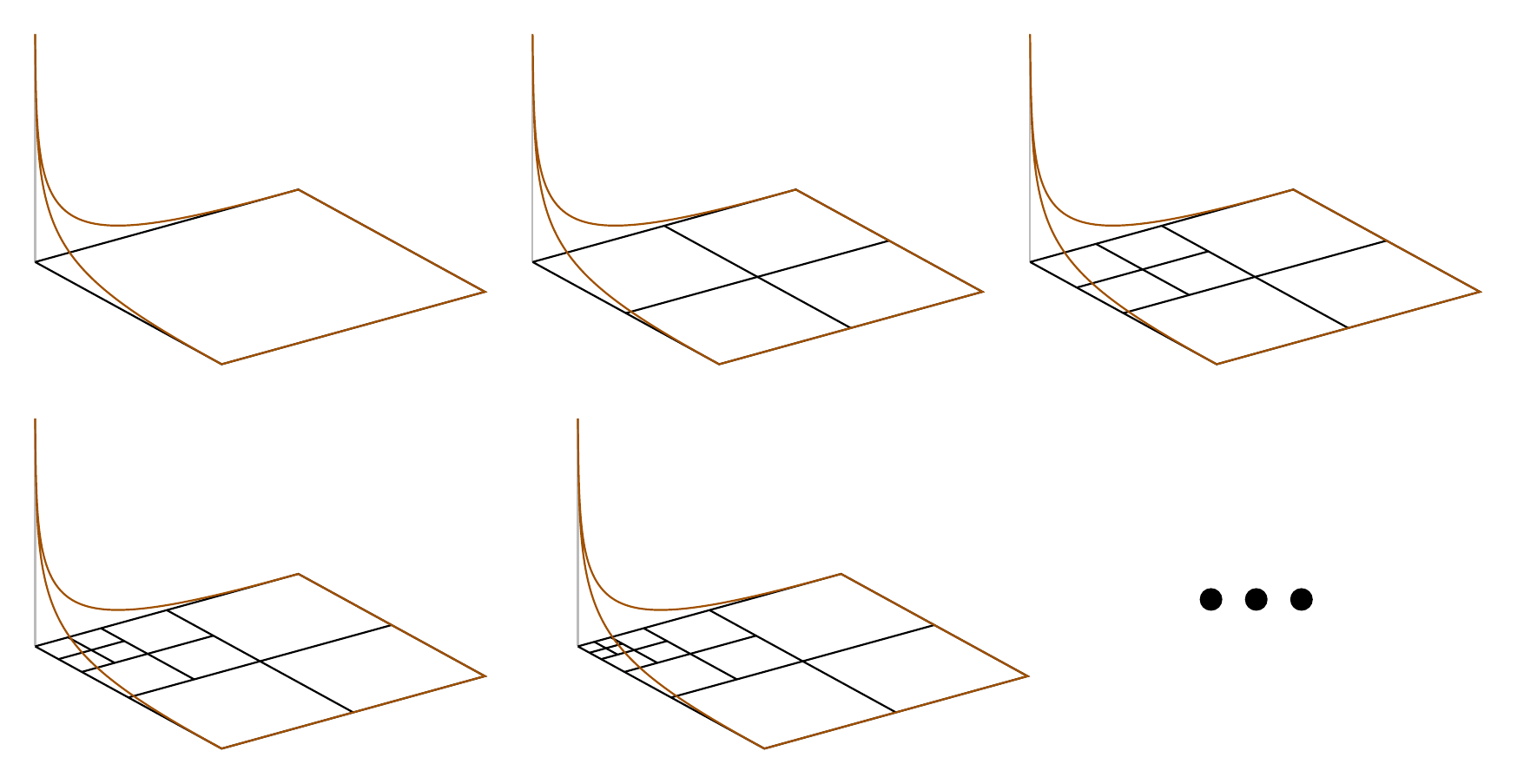}
    %\input{figures/ch2/rys2-1-example-mesh-constr.pdf_tex}
    %\includesvg[width=0.9\textwidth]{rys2-1-example-mesh-constr}
    \caption{Mesh refined over a point singularity \Aniarev{resulting from a gradient of} a two-dimensional function. The shape of the singularity $S$ is a point in the corner.}
    \label{ch2:example-refined}
\end{figure}

For this work, we employ $d$-dimensional hypercube elements (rectangles in 2D, hexahedrons in 3D, and octachorons in 4D), and we consider hierarchical polynomial basis functions of order $p$. Hierarchical basis functions are constructed by the tensor product of one-dimensional hierarchical shape functions glued together to obtain globally continuous piecewise polynomial functions with compact support.
We refer to \cite{Book2} for a more general description. \\
We consider shape functions over vertices, edges, faces, and interiors in three dimensions, while in four dimensions over vertices, edges, faces, hyperfaces, and interiors.\\
We identify nodes of the mesh with basis functions. We also consider nodes' support, defined as support of basis functions associated with the node. We have basis functions assigned to vertices, edges, faces, hyperfaces (in higher dimensions), and interiors. In general, \textcolor{black}{the support of the nodal function spans over all the adjacent elements sharing the node}. For example, in the case of two-dimensional regular mesh, the support of the vertex node spans into four elements sharing the node, the support of the edge node spans into two elements sharing the node, and the support of the interior node is equal to the single element.\\
We focus on $h$-adaptive meshes, where we employ the 1-irregularity rule, stating that an element can be broken only once without breaking its neighbors.
To construct an $h$-adaptive mesh around a singularity, we start with the one initial element and iteratively refine all elements that overlap with the singularity, ensuring that the $1$-irregularity rule is being followed. For example, Algorithm \ref{algo:mesh-construction} can be used to build $h$-adaptive mesh around a singularity with shape $S$ with refinement level $R$, where the symbol `` $[]$ '' stands for initiation as an empty array. \Aniarev{We assume that information on the
singularity location is available.}
An exemplary singularity mesh construction is shown in Figure \ref{ch2:example-refined}.

% \clearpage

\begin{algorithm}
    \begin{algorithmic}
    \Procedure{Constructmesh}{$RootElement, R, S$}
        \State $root \gets RootElement$
        \State $G \gets []$ \Comment{$G[r]$ contains all elements of refinement level $r$.}
        \State $G[0] \gets [root]$ \Comment{Initialize the mesh with single element.}
        \For{$r \gets 1$  \textbf{to} $R$} \Comment{\textbf{Step 1}: Refine all elements as necessary.}
            \State $H \gets G[r-1]$ \Comment{Only elements of refinement $r-1$ can be refined further.}
            \State $G[r-1] \gets []$
            \State $G[r] \gets []$
            \ForAll{$e \in H$}
                \If{$e \text{ overlaps } S$}
                    \If{$i > 1$}
                        \State $K \gets G[r-2]$
                        \State $G[r-2] \gets []$
                        \ForAll{$f \in K$}
                            \If{$f \text{ shares a vector with } e$}
                                \State $G[r-1] \gets \Call{concat}{G[r-1], \textsc{refine}(f)}$
                            \Else
                                \State $G[r-2] \gets \Call{append}{G[r-2], f}$
                            \EndIf
                        \EndFor
                    \EndIf
                    \State $G[r] \gets \Call{concat}{G[R], \textsc{refine}(e)}$ \Comment{If analyzed element overlaps with the singularity, refine it into smaller elements}

                \Else
                    \State $G[r-1] \gets \Call{append}{G[r-1], e}$ \Comment{...otherwise, leave it unrefined. }
                \EndIf
            \EndFor
        \EndFor
        \State $M \gets []$
        \For{$r \gets 0$ \textbf{ to } $R$}
            \State $M \gets \Call{concat}{M, G[r]}$
        \EndFor
        \State \Return $M$
    \EndProcedure
    \end{algorithmic}
    \caption{mesh construction over a singularity}
    \label{algo:mesh-construction}
\end{algorithm}

%\section{Sparse Gaussian elimination and matrix permutations based on element partition tree}

\subsection{Element partition tree}
\begin{figure}[!h]
\centering
\includegraphics[width=\textwidth]{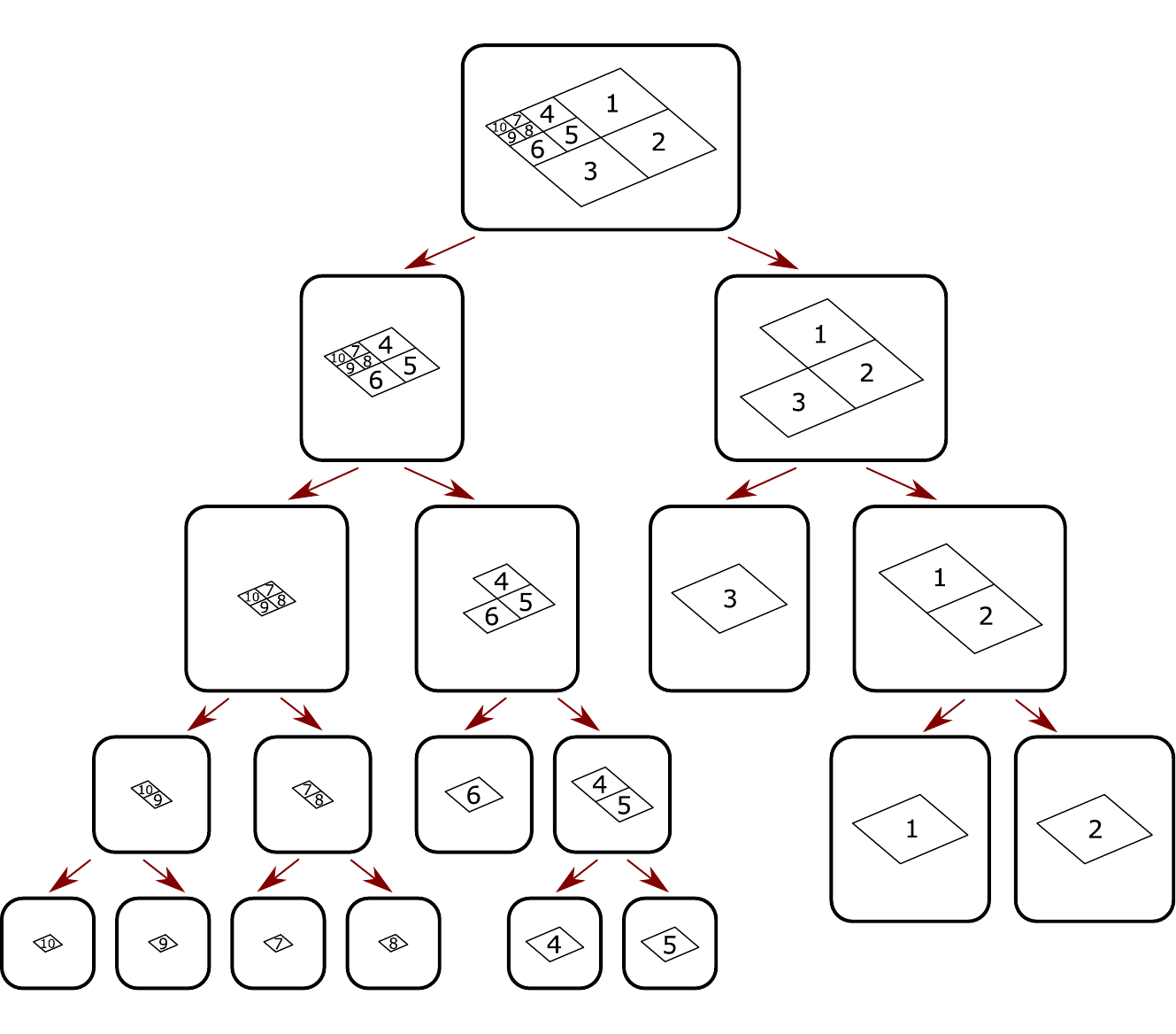}
\caption[Example of an element partition tree.]{An example for an element partition tree for a small adaptive mesh. The post-order traversal of the tree results in the ordering of elements $(10, 9, 7, 8, 6, 4, 5, 3, 1, 2)$.}
\label{ch2:element-partition-example}
\end{figure}

\Aniarev{The element partition tree is used for the construction of ordering for sparse matrix permutation in order to speed up the multi-frontal solver. Using some partitioning strategy, an element partition tree is created by recursive partitioning the mesh elements into two parts \cite {threeD1}, \cite{twoD1},\cite{twoD2}. An example of an element partition is presented in Figure \ref{ch2:element-partition-example}.
\\
An element partition tree for a mesh consisting of a set of elements $E$ is a binary tree defined as $T=(E, V,c_1,c_2,e)$ with the following properties:
\begin{enumerate}
\item $V$ is a set of tree nodes,
\item $c_1:V\rightarrow V$ and $c_2:V\rightarrow V$ are the functions assigning left and right child to a node, respectively,
\item $e:V\rightarrow P(A) $ is a function assigning subsets of a set of all elements of the mesh $E$ to a node,
\item the root called node $ROOT$ contains all elements of the mesh, in other words, $e(ROOT) = E$,
\item each node $node$ for which $|e(node)|=1$ is a leaf in the tree,
\item each node $node$ for which $|e(node)|>1$ has exactly two children; $c_1(n)$ and $c_2(n)$,
\item for each node $node$ $e(node) = e(c_1(node))\cup e(c_2(node))$ .
\end{enumerate}}

\subsection{Ordering generation and solving using generated ordering}
\begin{figure}[!h]
    \centering
    \def\svgwidth{\textwidth}
  \includegraphics[width=0.5\textwidth]{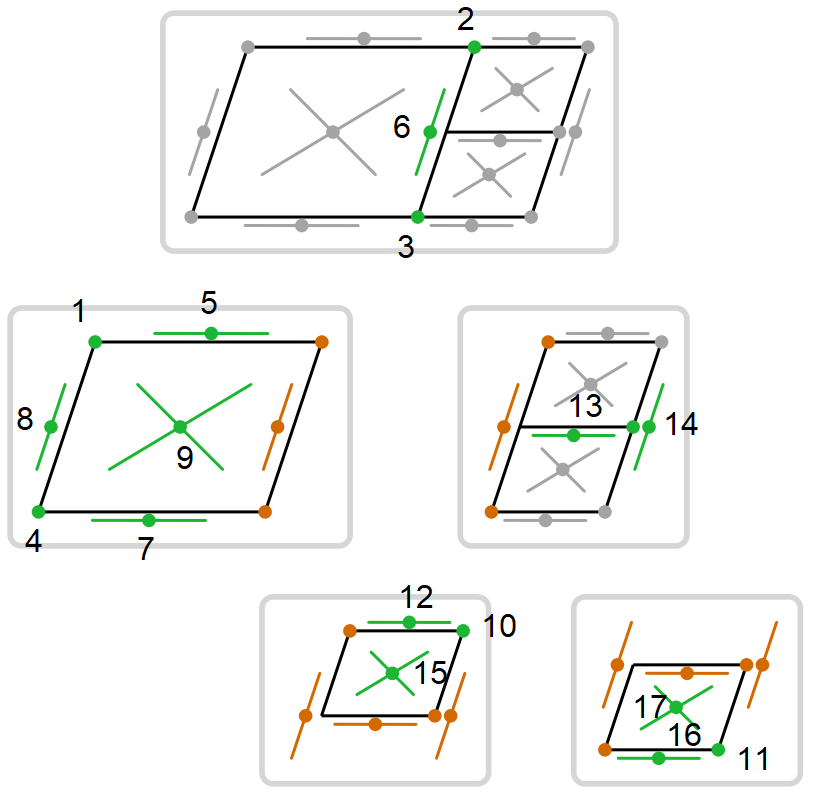}
    \caption[Element partition tree ordering generation.]{An example of element partition tree ordering generation. Nodes marked in green are eliminated at the given tree node. This example results in the permutation of $(1, 4, 5, 7, 8, 9, 10, 12, 15, 11, 16, 17, 13,14, 2, 3, 6)$. }
    \label{ch3:ordering-generation}
s\end{figure}

Given an element partition tree, we can generate a row elimination order for the matrix using a post-order traversal of the element partition tree. At each traversed partition tree node, we list all basis functions with support entirely contained by the tree node elements which have not been listed already. This produces a permutation (or ordering) of all nodes. \Aniarev{Figure \ref{ch3:ordering-generation} presents an exemplary element partition tree with denoted basis functions for each node, which generates the ordering $(1, 4, 5, 7, 8, 9, 10, 12, 15, 11, 16, 17, 13,14, 2, 3, 6)$.}

\subsection{Sparse Gaussian elimination}

We focus on the Gaussian elimination adapted to work over sparse matrices.  The algorithm solves a linear equation system represented by $Mx=A$ ($M$ is a sparse matrix and $A$ is a vector) by multiplying both sides by some ordering $P$ and solving $PMx = PA$ ($P$ being a permutation matrix representing some ordering).  In our case, the permutation of the matrix is based on the post-order traversal of the element partition tree.  The practical implementation of the algorithm is the multi-frontal solver \cite{DS1,DS2},  avoiding zeros in matrices by constructing the elimination tree internally based on the proposed ordering and the sparsity pattern of the matrix.  In our numerical experiments, we employ the multi-frontal solver from Octave.

\subsection{Time complexity of element partition tree based solvers}

This section shows that an element partition tree-based ordering gives a recursive formula for computation complexity that is easy to calculate. First, let us analyze the Gaussian elimination algorithm. It comprises two main steps--a first step corresponding to a forward elimination and a second step given by a backward substitution. The second step requires several operations proportional to the number of non-zero elements in the row form matrix. However, each non-zero element has to be non-zero at the beginning of the algorithm or originates from an operation performed in the first step. Therefore, the second step does not add anything to the computational complexity of the whole algorithm.
For practical reasons, most sparse matrix algorithms keep an element in the memory even if it has been modified to be $0$. For the sake of brevity,  we call a {\em non-zero element} to any element that is or has previously been set to a non-zero value, without regard to whether it is equal to $0$ at a given time.\\
The computational complexity of the first step can be analyzed as a sum of the complexities of eliminating rows for each element partition tree node. Let us make the following set of observations:
\begin{enumerate}
\item {A non-zero element in the initial matrix happens when the two basis functions corresponding to that row and column have overlapping supports. Let us call the graph created by considering the initial matrix to be an adjacency matrix of a graph as an {\em overlap graph}. Two graph nodes cannot be neighbors in an overlap graph unless the supports of their corresponding basis functions overlap. }
\item {When a row is eliminated,  the new non-zero elements are created on the intersection of columns and rows that has non-zero values in the eliminated row or corresponding column. If we analyze the matrix as a graph, then elimination of the row corresponding to a graph node produce edges between all pairs of nodes that were neighbors of the node being removed.}
\item {If at any given time during the forward elimination step a non-zero element exists on the intersection of a row and a column corresponding to two basis functions, then either those two basis functions have corresponding graph nodes that are neighbors in the overlap graph, or that there exists a path between those two nodes in the overlap graph that traverses only elements that have been eliminated already.}
\item {
All variables corresponding to the neighboring nodes of the graph node of a variable \Aniarev{$x_j$} in the overlap graph are either:
\begin{enumerate}
\item {listed in one of the element partition tree nodes that are descendants of the element partition tree node listing the variable \Aniarev{$x_j$} – and those variables are eliminated already by the time this variable is eliminated, or}
\item {listed in the same element partition tree node as the variable \Aniarev{$x_j$}, or}
\item { having the support of the corresponding basis function intersected by the boundary of the submesh of the element partition tree node containing the variable \Aniarev{$x_j$} – those graph nodes are listed in one of the ancestors of the element partition tree node listing the variable \Aniarev{$x_j$}. }
\end{enumerate}
Thus, in the overlap graph, there are no edges between nodes that belong to two different element partition tree nodes that are not in an ancestor-descendant relationship. At the same time, any path that connects a pair of non-neighboring nodes in the overlap graph has to go through at least one graph node corresponding to a variable that is listed in a common ancestor of the element partition tree nodes containing the variables from that pair of nodes. }
\end{enumerate}
% \item {All neighbors of a mesh node in an overlap graph will either already be eliminated by the time that node is eliminated, or are to be listed while traversing the same tree node, or are on the interface of that tree node and will be eliminated later. This means that there will not be any non-zero values on the intersection of the removed row/column and columns/rows corresponding to mesh nodes that do not overlap that tree node's submesh.}
These observations lead to a conclusion that a removal of a single row requires no more than \Aniarev{$(z-1)z $} subtractions, where \Aniarev{$z$ is the amount of variables related with graph nodes that were not removed earlier, have overlapping support and are listed at the same tree node.} Last leads to the following recursive formula for the complexity of removing all mesh nodes at a single tree node:
\Aniarev{
\begin{equation}
\begin{array}{r @{\;} l}
T(node) &= \displaystyle T(child_1(node)) + T(child_2(node)) + \sum_{i=1}^{a}(b-i)(b-i+1) \smallskip \\
&= \displaystyle T(child_1(node)) + T(child_2(node)) + \frac{1}{3}a(a^2 + 6 b + 3 b^2 + 3 a b+ 3 a + 2) \smallskip \\
&= \displaystyle T(child_1(node)) + T(child_2(node)) + \mathcal{O}(ab^2),
\end{array}
\end{equation}}
where
\Aniarev{$a$ denotes the number of variables removed at a given tree node, and
$b$ is equal to $a$ plus the number of variables on the interface of that tree node. \\
In the following sections,  we denote the complexity of the removal of $a$ variables belonging to a tree node from $b$ variables }having overlapping support with the corresponding node as:
\begin{equation}
    C_r(a,b) = \frac{1}{3}a(a^2 + 6 b + 3 b^2 + 3 a b+ 3 a + 2) = \mathcal{O}(ab^2).
\end{equation}
\section{Computational complexities of space--time formulations and time--marching schemes}
\begin{figure}[ht]
  \includegraphics[width=0.75\textwidth]{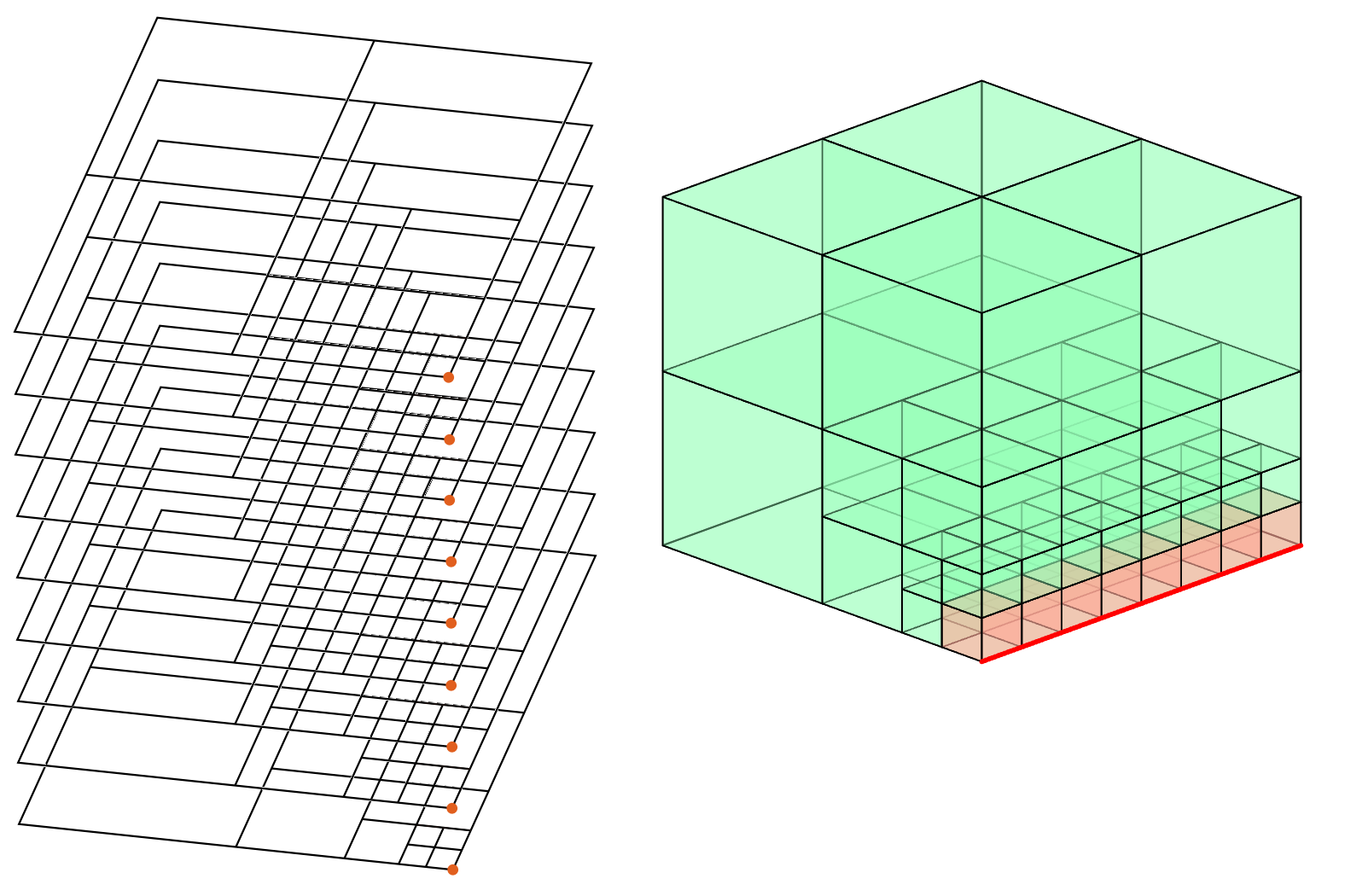}
  \centering
  \caption{The refinement towards a space--time ``edge'' resulting from a point traveling through space and time.  This space--time refinement pattern corresponds to a sequence of spatial meshes refined towards ``points'' located on the space--time edge at particular time moments. }
  \label{ch7:edgepoint}
\end{figure}

This section considers some possible space--time mesh refinement patterns resulting from different space--time singularities (cf.~\ref{ApA} and \ref{ApB}). We look at the structure of the $d$-dimensional space--time refined mesh. Namely, we focus on
\begin{itemize}
\item space--time ``edge'',  resulting from a point traveling through space and time,
\item space--time ``face'',  resulting from an edge traveling through space and time,
\item space--time ``hyperface'',  resulting from a face traveling through space and time.
\end{itemize}
We consider these space--time refined meshed in three-dimensions ($d=3$) and in four-dimensions ($d=4$).
These space--time refined $d$-dimensional meshes correspond to the following
sequences of $d-1$ dimensional refined meshes employed by the time--marching scheme:
\begin{itemize}
\item sequence of spatial meshes refined towards ``points'' located on the space--time edge at particular time moments,
\item sequence of spatial meshes refined towards ``edges'' located on the space--time face at a particular time moment,
\item sequence of spatial meshes refined towards ``faces'' located on the space--time hyperface at a particular time moment.
\end{itemize}
This correspondence is illustrated in Figure \ref{ch7:edgepoint} for the three-dimensional space--time mesh with edge singularity, and the resulting sequence of two-dimensional meshes refined to the corresponding points.
We assume that the time step size in the time--marching scheme is equal to the diameter along the time axis of the smallest element in the space--time mesh.
\subsection{Computational complexities of direct solvers for space--time formulations and time--marching schemes}
We compare the computational complexity of the direct solver executed over the $d$ dimensional space--time domain to the computational complexity of the direct solver executed multiple times over $d$-1 dimensional meshes within the time--marching scheme.
In Table ~\ref{ch5:SpaceTimeComplexitiesTable} \Aniarev{ the mesh dimensions and direct solver time complexity for space--time adaptive meshes are presented}. We include the 4D space--time uniform mesh, 4D space--time mesh refined towards the hyperface singularity, 4D space--time mesh refined towards the face singularity, 4D space--time mesh refined towards edge singularity, as well as 3D space--time uniform mesh, 3D space--time mesh refined towards face singularity, and 3D space--time mesh refined towards edge singularity. We compute \Aniarev{sizes} $N$ of all these grids and express it by the number of refinement levels $r$. Additionally, we provide the estimates for the direct solver execution time. \Aniarev{The estimations presented in Table ~\ref{ch5:SpaceTimeComplexitiesTable} are based on theoretical results presented in Appendix A and Appendix B, see Table ~\ref{ch5:table1} and  ~\ref{ch5:complexitiesTable}}.
An interesting observation is that the computational complexity of the direct solver executed for the space--time mesh refined towards an edge is linear ${\cal O}(N)$.\\
Next, in Table \ref{ch5:TimeMarchingComplexitiesTable} we evaluate the computational complexity for time--marching scheme. For each $d$-dimensional space--time refined toward the $q$-singularity, we construct a sequence of $2^r$ meshes. Each of the grids from the sequence is $d-1$ dimensional, refined towards $q-1$ singularity. Such a sequence of grids provides the solution with similar accuracy to the space--time grid. We estimate the dimensions $n$ of the spatial meshes from the sequence and the computational complexity of the direct solver executing $2^r$ times, once for each grid from the sequence. \Aniarev{The estimations presented in Table \ref{ch5:TimeMarchingComplexitiesTable} are based on theoretical results presented in Appendix A and Appendix B, see Table ~\ref{ch5:table1} and  ~\ref{ch5:complexitiesTable}}.
An interesting observation is that the computational complexity of the time--marching scheme corresponding to the space--time mesh refined towards an edge is $r$ times higher than the space--time complexity. However, the time marching scheme is cheaper for all other space--time grids.
\begin{table}[h!]
\centering
\begin{tabular}{ |c|c|c|c| }
%& \multicolumn{3}{|c|}{Singularity type} \\
 \hline
space--time mesh  & space--time mesh size        & space--time mesh \\
 & $N$ & direct solver complexity  \\
\hline
4D uniform & ${\cal O}(2^{4r})$ &   \inred{${\cal O}(N^{\frac{9}{4}})={\cal O}(2^{9r})$} \\
\hline
4D hyperface & ${\cal O}(8^r)={\cal O}(2^{3r})$ & \inred{${\cal O}(N^2)={\cal O}(2^{6r})$} \\
\hline
4D face & ${\cal O}(4^r)={\cal O}(2^{2r})$  & \inred{${\cal O}(N^{\frac{3}{2}})={\cal O}(2^{3r})$ }\\
\hline
4D edge & ${\cal O}(2^r)$ & \inred{${\cal O}(N)={\cal O}(2^{r})$} \\
\hline
3D uniform & ${\cal O}(2^{3r})$ &   \inred{${\cal O}(N^{2})={\cal O}(2^{6r})$} \\
\hline
3D face & ${\cal O}(4^r)={\cal O}(2^{2r})$ & \inred{${\cal O}(N^{\frac{3}{2}})={\cal O}(2^{3r})$} \\
\hline
3D edge & ${\cal O}(2^r)$ & \inred{${\cal O}(N)={\cal O}(2^{r})$} \\
\hline
\end{tabular}
\caption{space--time formulations, mesh dimensions and direct solver complexities}
\label{ch5:SpaceTimeComplexitiesTable}
\end{table}
\begin{table}[h!]
\centering
\begin{tabular}{ |c|c|c|c| }
%& \multicolumn{3}{|c|}{Singularity type} \\
 \hline
space--time mesh & Sequence of & Single spatial & Total direct solver \\
 & spatial meshes & mesh size & complexity for a sequence \\
\hline
4D uniform
& $\inred{2^r \times} $ 3D uniform & $n={\cal O}(2^{3r})$ & \inred{$2^r \times {\cal O}(n^{2})=$} \\
\inblue{$N={\cal O}(2^{4r})$} & & & \inblue{${\cal O}(2^{7r})$} \\
\hline
4D hyperface
& $\inred{2^r \times} $ 3D face & $n={\cal O}(2^{2r})$
& \inred{$2^r \times{\cal O}( n^{\frac{3}{2}})=$}
\\
\inblue{$N={\cal O}(2^{3r})$}  & & & \inblue{${\cal O}(2^{4r})$} \\
\hline
4D face
& $\inred{2^r \times} $ 3D edge & $n={\cal O}(2^{r})$  &
\inred{$2^r \times{\cal O}( n)=$}
\\
\inblue{$N={\cal O}(2^{2r})$}  &   & & \inblue{${\cal O}(2^{2r})$} \\
\hline
4D edge
& $\inred{2^r \times} $ 3D point & $n={\cal O}(r)$ &
\inred{$2^r \times{\cal O}(r)=$}
\\
\inblue{$N={\cal O}(2^{r})$}    &  & & \inblue{${\cal O}(r2^r)$} \\
\hline
3D uniform
& $\inred{2^r \times} $ 2D uniform & $n={\cal O}(2^{2r})$ & \inred{$2^r \times {\cal O}(n^{\frac{3}{2}})=$} \\
\inblue{$N={\cal O}(2^{3r})$}   &  & & \inblue{${\cal O}(2^{4r})$} \\
\hline
3D face
& $\inred{2^r \times} $ 2D edge & $n={\cal O}(2^{r})$ &
\inred{$2^r \times{\cal O}(2^r)=$}
 \\
\inblue{$N={\cal O}(2^{2r})$}    &  & & \inblue{${\cal O}(2^{2r})$} \\
\hline
3D edge
& $\inred{2^r \times} $ 2D point & $n={\cal O}(r)$  &
\inred{$2^r \times{\cal O}(r)=$}
\\
\inblue{$N={\cal O}(2^{r})$}    & & & \inblue{${\cal O}(r2^r)$} \\
\hline
\end{tabular}
\caption{time--marching schemes corresponding to space--time formulations}
\label{ch5:TimeMarchingComplexitiesTable}
\end{table}

\subsection{Computational complexities of iterative solvers for space--time formulations and time--marching schemes}
Finally, we focus on the iterative solver. The computational complexities of the iterative solver for both space--time grids and the time--marching scheme are estimated as in Table \ref{ch5:IterativeSolvers}. The computational complexity for the space--time grid is equal to ${\cal O}(N_{iter}\times N)$ where $N$ is the size of the space--time grid, and $N_{iter}$ is the number of iterations of the iterative solver for the space--time grid.
The computational complexity for the space--time grid is equal to ${\cal O}(2^r \times n_{iter}\times n)$ where $2^r$ is the number of steps of the time--marching scheme, $n$ is the size of the spatial grid from the sequence, and $n_{iter}$ is the number of iterations on the spatial grid.
It is obvious to assume that $N_{iter} >> n_{iter}$, but the exact numbers of iterations are problem-dependent.
The computational complexity of the iterative solver for space--time grids is higher than the computational complexity of the iterative solver for the time--marching scheme.
An interesting case is the space--time grid refined to the space--time edge, following the trajectory of the point object. In this case, the complexity of iterative solver for space--time grid is ${\cal O}(N_{iter}\times N)={\cal O}(N_{iter} \times 2^r)$ while the complexity of the iterative solver for time marching grids is ${\cal O}(2^r \times n_{iter}\times n)={\cal O}(n_{iter} \times r2^r)$. The iterative solver for space--time grid is cheaper than for the time marching scheme, if $N_{iter} < r \times n_{iter}$.
\begin{table}[h!]
\centering
\begin{tabular}{ |c|c|c|c| }
%& \multicolumn{3}{|c|}{Singularity type} \\
 \hline
space--time  & Sequence of & space--time & time--marching scheme \\
mesh & spatial meshes & iterative & iterative \\
& & solver complexity & solver complexity \\
\hline
4D uniform  & $\inred{2^r \times} $ 3D uniform &
& \inred{${\cal O}(2^{r}\times n_{iter} \times 2^{3r})=$} \\
\inblue{$N={\cal O}(2^{4r})$} & &\inblue{${\cal O}(N_{iter} \times 2^{4r})$} & \inblue{${\cal O}(n_{iter} \times 2^{4r})$} \\
\hline
4D hyperface  & $\inred{2^r \times} $ 3D face &
& \inred{${\cal O}(2^{r}\times n_{iter} \times 2^{2r})=$}
\\
\inblue{$N={\cal O}(2^{3r})$} & & \inblue{${\cal O}(N_{iter} \times 2^{3r})$} & \inblue{${\cal O}(n_{iter} \times 2^{3r})$} \\
\hline
4D face   & $\inred{2^r \times} $ 3D edge &   &
\inred{${\cal O}(2^r \times n_{iter} \times 2^{r})=$}
\\
\inblue{$N={\cal O}(2^{2r})$} & & \inblue{${\cal O}(N_{iter} \times 2^{2r})$} & \inblue{${\cal O}( n_{iter} \times 2^{2r})$}\\
\hline
4D edge  & $\inred{2^r \times} $ 3D point &  &
\inred{${\cal O}(2^r \times n_{iter} \times r) = $}
\\
\inblue{$N={\cal O}(2^{r})$} & & \inblue{${\cal O}(N_{iter}\times 2^r )$} & \inblue{${\cal O} (n_{iter} \times r 2^r )$} \\
\hline
3D uniform  & $\inred{2^r \times} $ 2D uniform &
& \inred{${\cal O}(2^{r}\times n_{iter} \times 2^{2r})=$} \\
\inblue{$N={\cal O}(2^{3r})$} & &  \inblue{${\cal O}(N_{iter} \times 2^{3r})$}  & \inblue{${\cal O}(n_{iter} \times 2^{3r})$} \\
 \hline
3D face   & $\inred{2^r \times} $ 2D edge &  &
\inred{${\cal O}(2^r \times n_{iter} \times 2^{r})=$}
 \\
\inblue{$N={\cal O}(2^{2r})$} & &\inblue{${\cal O}(N_{iter} \times 2^{2r})$} & \inblue{${\cal O}( n_{iter} \times 2^{2r})$}\\
\hline
3D edge   & $\inred{2^r \times} $ 2D point &  &
\inred{${\cal O}(2^r \times n_{iter} \times r) = $}
\\
\inblue{$N={\cal O}(2^{r})$} & &  \inblue{${\cal O}(N_{iter}\times 2^r )$} & \inblue{${\cal O} (n_{iter} \times r 2^r )$} \\
\hline
\end{tabular}
\caption{Comparison of complexities of iterative solvers for space--time formulations and time--marching schemes.}
\label{ch5:IterativeSolvers}
\end{table}
Additionally, we compare the computational complexity of the iterative solver executed over the $d$ dimensional space--time domain to the computational complexity of the iterative solver executed multiple times over the $d$-1 dimensional meshes during the time--marching scheme.
\subsection{Impact of polynomial order of approximation}
Notice that we have not included the $p$-factor in the computational complexity estimates for the arbitrary shape of singularity. However, we can easily estimate the computational complexity of the static condensation performed at the beginning of the computations with a higher-order finite element method with hierarchical basis functions.\\
The computational complexity of the static condensation over a single element is equal to ${\cal O}(p^{3d})$. We are eliminating the degrees of freedom from the interior of the element; we have $(p-1)^d$ degrees of freedom there; the matrix is dense, so the complexity of elimination is the cube of the number of degrees of freedom. In other words, the complexity of static condensation for a grid with $N$ elements is ${\cal O}(Np^{3d})$. With this observation in mind, we can estimate the complexities of static condensations for space--time and time--marching grids. It is illustrated in Table \ref{ch5:staticcondensations}. The complexity of static condensations is higher for the space--time mesh.
\revision{However, for the space--time formulation, we can always execute all the elemental computations at the same time, in parallel, while for the time--marching scheme, we can only run the static condensation over a single time-step mesh in parallel. In this sense, the space--time method can outperform the time--marching scheme.}
\begin{table}[h!]
\centering
\begin{tabular}{ |c|c|c|c| }
%& \multicolumn{3}{|c|}{Singularity type} \\
 \hline
space--time & Sequence of & space--time & time--marching scheme \\
mesh  & spatial meshes & complexity of static & complexity of static\\
& & condensation & condensation \\
\hline
4D uniform  & \inred{$2^r$ steps $\times $} & & \inred{($2^r$ steps)} \\
\inblue{$N={\cal O}(2^{4r})$} & 3D uniform  & \inblue{${\cal O}(p^{12}2^{4r})$}
& \inblue{${\cal O}(\inred{2^{r}} \times p^9 2^{3r})$} \\
\hline
4D hyperface  &  \inred{$2^r$ steps $\times $} & & \inred{($2^r$ steps)}
\\
\inblue{$N={\cal O}(2^{3r})$} &  3D face & \inblue{${\cal O}(p^{12}2^{3r})$}
& \inblue{${\cal O}(\inred{2^{r}} \times p^9 2^{2r})$} \\
\hline
4D face   &  \inred{$2^r$ steps $\times $}& & \inred{($2^r$ steps)}
\\
\inblue{$N={\cal O}(2^{2r})$} &  3D edge & \inblue{${\cal O}(p^{12} 2^{2r})$}  &
\inblue{${\cal O}(\inred{2^{r}} \times p^{9} 2^{r})$} \\
\hline
4D edge  &  \inred{$2^r$ steps $\times $} & & \inred{($2^r$ steps)}
\\
\inblue{$N={\cal O}(2^{r})$} &  3D point & \inblue{${\cal O}(p^{12}  2^r )$} &
\inblue{${\cal O}(\inred{2^{r}}  \times p^9r) $} \\
\hline
3D uniform  & \inred{$2^r$ steps $\times $} & & \inred{($2^r$ steps)} \\
\inblue{$N={\cal O}(2^{3r})$} &  2D uniform & \inblue{${\cal O}(p^9 2^{3r})$}
& \inblue{${\cal O}(\inred{2^{r}} \times p^6 2^{2r})$} \\
 \hline
3D face   & \inred{$2^r$ steps $\times $} & & \inred{($2^r$ steps)}
 \\
\inblue{$N={\cal O}(2^{2r})$} &   2D edge & \inblue{${\cal O}(p^9 2^{2r})$} &
\inblue{${\cal O}(\inred{2^{r}} \times p^6 2^{r})$} \\
\hline
3D edge   & \inred{$2^r$ steps $\times $} & & \inred{($2^r$ steps)}
\\
\inblue{$N={\cal O}(2^{r})$} & 2D point  &  \inblue{${\cal O}(p^9 2^r )$} &
\inblue{${\cal O}(\inred{2^r} \times p^6 r)$} \\
\hline
\end{tabular}
\caption{Comparison of complexities of static condensations for space--time formulations and time--marching schemes.}
\label{ch5:staticcondensations}
\end{table}
\section{Numerical results}
In this section we provide numerical experiments for verification of the theoretical findings summarized in Tables \ref{ch5:SpaceTimeComplexitiesTable}-\ref{ch5:IterativeSolvers}.
We employ our Octave codes to generate the structure of $d$-dimensional computational grids with $q$-dimensional singularity, generated using $r$ refinement levels.
For simplicity of implementation, we generate our matrices by looking at relations between finite elements. Rows and columns in matrices correspond to finite elements. Non-zero entries in a row mean that two elements, one related to the row, and one related to the column, are adjacent through a $d-1$ dimensional face.
This way of generating matrices influences the computational complexity constant, ignoring the polynomial order of approximation, but the dependence on $N$ equal here the number of elements is of the same order as if we include all the relations of the basis functions.\\
We run experiments using Octave on a Linux cluster node equipped with 2.4GHz processor with 64 GB of RAM. We cannot factorize more than eight refinements for the face singularity in four dimensions and nine refinement levels for the face singularity in three dimensions because of a lack of memory during the factorization process.
The comparisons of execution times for four-dimensional face singularity versus a sequence of three-dimensional edge singularities are presented in Table \ref{tab:4Dface}.
The comparisons of execution times for three-dimensional face singularity versus a sequence of two-dimensional edge singularities are presented in Table \ref{tab:3Dface}.
We employ AMD ordering and multi-frontal solver as implemented in the Octave, e.g., for the face singularity with nine refinement levels in three dimensions, we run our matrix generation script:\\

\texttt{F9\_3=Face(9,3);}

\texttt{N = 87381}

We compute the AMD permutation

\texttt{p=amd(F9\_3);}

and we plug it into the LU factorization, measuring the execution time

\texttt{tic; lu(F9\_3(p,p)); toc}

\texttt{Elapsed time is 46.2016 seconds.}\\

\noindent
We know this ordering is different from the one proposed in our paper. Nevertheless, the results show up to one order of magnitude times faster execution times of time--marching schemes with edge singularities than one call for the space--time domain with the face singularity.\\
The comparisons of execution times for four-dimensional edge singularity versus a sequence of three-dimensional point singularities are presented in Table \ref{tab:4Dedge}.
We can perform 15 refinements over the space--time mesh this time with the Octave solver. The space--time mesh with a multi-frontal solver is, in this case, faster than the time marching scheme, up to the 13 refinement level. With 14 or 15 refinements, the complexity of processing the space--time mesh is higher than the complexity of processing the time--marching solver.
%
%
%\begin{figure}[ht]
%  \includegraphics[width=0.33\textwidth]{Point4_2.png}\includegraphics[width=0.33\textwidth]{Point4_3.png}\includegraphics[width=0.33\textwidth]{Point4_4.png}
%  \centering
%  \caption{Structure of the matrix for two-, three-, and four- dimensional meshes refined towards point singularity with $r=4$.}
%  \label{ch7:point4d}
%\end{figure}

\begin{table}[h!]
\centering
\begin{tabular}{ |c|c|c|c|c|c|}
 \hline
space--time & r & N & space--time & Sequence of & time--marching scheme \\
mesh & & & solver time [s] & spatial meshes & solver time [s] \\
\hline
4D face & 6 & 4095 & 0.56 & $64\times$ 3D edge &
$64 \times 0.0037 =0.0236$ \\
\hline
4D face & 7 & 16383 & 13.79 & $128\times $ 3D edge &
$128 \times 0.024 =3.072$ \\
\hline
4D face & 8 & 65535 & 234 & $256\times $ 3D edge &
$256 \times 0.063 =16.128 $ \\
\hline
\end{tabular}
\caption{Execution times for four dimensional mesh refined towards face singularity and the corresponding sequence of three-dimensional meshes refined towards edge singularities. }.
\label{tab:4Dface}
\end{table}

\begin{table}[h!]
\centering
\begin{tabular}{ |c|c|c|c|c|c|}
 \hline
space--time & r & N & space--time & Sequence of & time--marching scheme \\
mesh & & & solver time [s] & spatial meshes & solver time [s] \\
\hline
3D face & 6 & 1365 & 0.044 & $64\times $ 2D edge &
$64 \times 0.0008=0.05$ \\
\hline
3D face & 7 & 5461 & 0.28 & $128\times $ 2D edge &
$128 \times 0.001=0.128$ \\
\hline
3D face & 8 & 21845 & 2.55 & $256\times $ 2D edge &
$256 \times 0.003=0.76$ \\
\hline
3D face & 9 & 87381 & 46.20 & $512\times $ 2D edge &
$512 \times 0.0074=3.78$ \\
\hline
\end{tabular}
\caption{Execution times for three dimensional mesh refined towards face singularity and the corresponding sequence of two-dimensional meshes refined towards edge singularities. }
\label{tab:3Dface}
\end{table}
%\begin{figure}[ht]
%  \includegraphics[width=0.33\textwidth]{Edge3_2.png}\includegraphics[width=0.33\textwidth]{Edge3_3.png}\includegraphics[width=0.33\textwidth]{Edge3_4.png}
%  \centering
%  \caption{Structure of the matrix for two-, three-, and four- dimensional meshes refined towards edge singularity with $r=3$.}
%  \label{ch7:edge4d}
%\end{figure}
%
%
%
%\begin{figure}[ht]
%  \includegraphics[width=0.33\textwidth]{Face3_3.png}\includegraphics[width=0.33\textwidth]{Face3_4.png}
%  \centering
%  \caption{Structure of the matrix for three-, and four- dimensional mesh refined towards face singularity with $r=3$.}
%  \label{ch7:face4d}
%\end{figure}
\begin{table}[h!]
\centering
\begin{tabular}{ |c|c|c|c|c|c|}
 \hline
space--time & r & N & space--time & Sequence of & time--marching scheme \\
mesh & & & solver time [s] & spatial meshes & solver time [s] \\
\hline
4D edge & 6 & 434 & 0.008 & $64\times$ 3D point &
$64 \times 0.0004 =0.0256$ \\
\hline
4D edge & 7 & 882 & 0.026 & $128\times $ 3D point &
$128 \times 0004=0.0512$ \\
\hline
4D edge & 8 & 1778 & 0.069 & $256\times $ 3D point &
$256 \times 0.0005 =0.128 $ \\
\hline
4D edge & 9 & 3570 & 0.15 & $512\times $ 3D point &
$512 \times 0.0006 =0.3 $ \\
\hline
4D edge & 10 & 7154 & 0.38 & $1024\times $ 3D point &
$1024 \times 0.0006 =0.61 $ \\
\hline
4D edge & 11 & 14332 & 1.01 & $2048\times $ 3D point &
$2048\times 0.0007 =1.43 $ \\
\hline
4D edge & 12 & 28665 & 2.45 & $4096\times $ 3D point&
$4096\times 0.0008 =3.27 $ \\
\hline
4D edge & 13 & 57330 & 5.52 & $8192\times $ 3D point&
$8192\times 0.0007 =5.73 $ \\
\hline
4D edge & 14 & 114674 & 13.59 & $16384\times $ 3D point&
$16384\times 0.0007 =11.46 $ \\
\hline
4D edge & 15 & 229362 & 38.96 & $32768\times $ 3D point&  $32768\times 0.0009 =29.49 $ \\
%\hline
%4D edge & 16 & 458738 & 120 & $65536\times $ 3D point&  $65536\times 0.001 =65.536 $ \\
\hline
\end{tabular}
\caption{Execution times for four dimensional mesh refined towards edge singularity and the corresponding sequence of three-dimensional meshes refined towards point singularities. }.
\label{tab:4Dedge}
\end{table}
\section{Conclusions}
To estimate the computational complexity of full space--time formulations and time--marching schemes for hypercubic elements, we simulate several possible scenarios for the resulting matrices when considering adaptivity toward singularities.
In particular, we consider refinements towards the point, edge, face, and hyperface singularities over space--time mesh. In our idealized case, we refine all the elements that contain the prescribed point, edge, face or hyperface singularity.
Thus, we obtain several representative refined $d$-dimensional computational meshes, where we assume that we perform refinements towards $q$-dimensional manifold ($q<d$) representing the singularities.
For each of these representative meshes, we estimate the number of degrees of freedom (second column in Table 1), the computational complexity of the multi-frontal solver (third column in Table 1), the computational complexity of the iterative solver (third column in Table 3), and the computational complexity of the static condensation (third column in Table 4).
On the other hand, we generated a sequence of refined $d-1$ dimensional computational meshes, representing the ``cross-sections" of the $d$ dimensional space--time mesh.
In this case, we performed refinements towards $q-1$-dimensional manifold representing the cross-section of the $q$-dimensional singularity.
Our theoretical estimations and the numerical experiments imply that the time--marching scheme is competitive only in this idealized case when the number of time steps is equal to the element size in the time dimension.
This concerns the computational complexity of the multi-frontal solver (third column in Table 1 versus the fourth column in Table 2), the iterative solver (third column in Table 3 and fourth column in Table 3), and the static condensation (third and fourth column in Table 4).
We also present numerical experiments, confirming the predicted theoretical behaviors.

We understand that our assumptions are the best possible idealistic scenarios.
In the real life applications, the computational complexity of space--time formulation is more competitive due to
\begin{itemize}
\item Increased number of time steps in higher-order and accuracy time marching schemes, where the time-step size can be actually smaller than the temporal dimension of the smallest elements in the space--time mesh.
\item Extensive parallelization of the computational process, where for example the static condensation for the space--time formulation can be performed fully in parallel, and the static condensation of the time--marching scheme has the limitation of the single size of the time-step mesh.
\item Sequential nature of the time--marching scheme, where the iterative solver has to be executed in a sequence for each time-step mesh, and it cannot be parallelized once for the entire computational space--time mesh.
\item The cost of generation of the refined computational meshes is in general ignored in our estimations (it is assumed to be linear), while in general, it is an iterative procedure that requires several solves, and in the space--time setup it can be performed once for the entire mesh, but in the time--marching scheme it has to be performed for each time-step mesh.
\end{itemize}
Nevertheless, our estimate constitutes the lower bounding case of the computational complexities for both space--time and time--marching schemes.

%As expected, space--time formulations performed on adaptive grids in the simulated scenarios usually are more expensive than old-fashioned time marching schemes for both direct and iterative solvers. We confirmed the last affirmation theoretically and experimentally. We summarized our findings concerning direct solvers in Tables 3-4, where we show the time complexity of the direct solvers on space--time grids (Table 3) and the time complexities of corresponding time marching schemes (Table 4). We also summarized our findings concerning iterative solvers in Table 5. The only exception is when we model a point object traveling in space (see the fourth and seventh row in Table 4 corresponding to the space--time edge singularity and sequence of spatial point singularities). We assume that the number of iterations of the iterative solver grows with problem size. Namely, $n_{iter} << N_{iter}$ where the first one corresponds to the number of iterations of an iterative solver run on the spatial mesh from a sequence of meshes, and the former one corresponds to the number of iterations of an iterative solver run on the space--time mesh. This time, space--time formulations are also more expensive, unless for the traveling point objects, where we have $N_{iter} < r \times n_{iter}$ where $r$ is the refinement level.
%In Table 6, we present computational complexities of static condensations, often performed before calling the solver \cite{static}. In this case, the complexities of static condensation of space--time formulation are always higher than for the time--marching scheme.
In our future work, we plan to extend our computational complexity estimates into parallel distributed memory \cite{wozniak1} and shared-memory machines \cite{wozniak2}, considering direct and iterative solvers. Our future work will also involve the computational complexity analysis of the hierarchical matrices solvers (H-matrices) in the space--time setup. In general, the complexity of processing these matrices is proportional to the problem size $N$, multiplied by an additional factor that grows with the dimension size.
\section*{Acknowledgments}
This work from the European Union's Horizon 2020 research and innovation programme under the Marie Sklodowska-Curie grant agreement No 777778 (MATHROCKS).
The work of SR has also been supported by the Chilean grant ANID Fondecyt No 3210009.
\appendix

\section{Computational complexity for h-adapted meshes towards a point singularity}\label{ApA}
This Appendix analyzes the time complexity of the direct solver being run on hierarchical meshes adapted toward a point singularity.
\subsection{\Aniarev{h-adapted mesh towards a point singularity }}\label{ch4:sec:def}
\begin{figure}[h!]
  \centering
  \caption{Examples of point singularities.}
  \begin{subfigure}[b]{0.9\textwidth}
  \includegraphics[width=\linewidth]{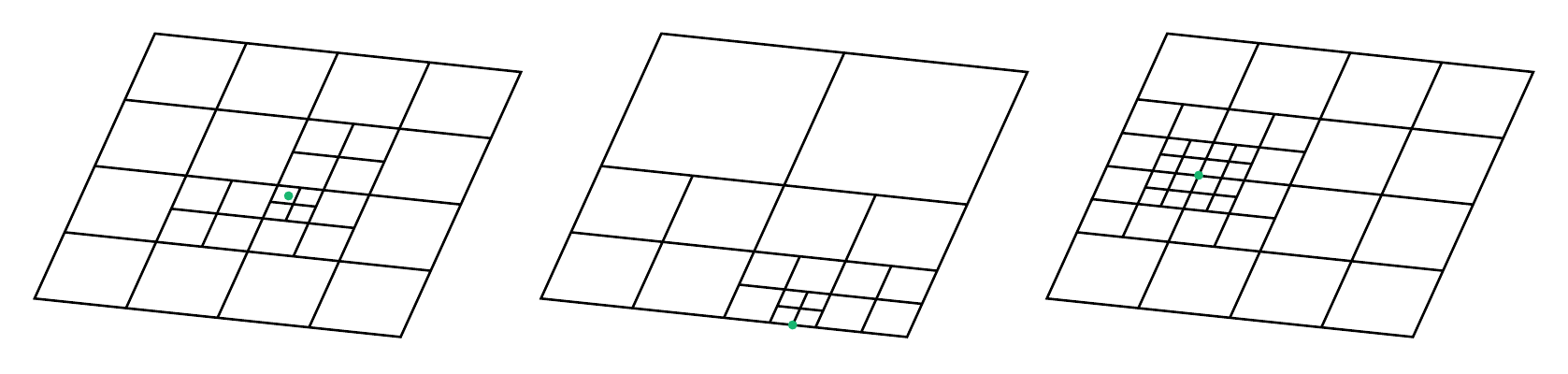}
       \caption{Examples of $h$-adaptive singularity meshes in 2D. Green dots denote the singularity position.}
    \end{subfigure}
    \begin{subfigure}[b]{0.6\textwidth}
    \includegraphics[width=\linewidth]{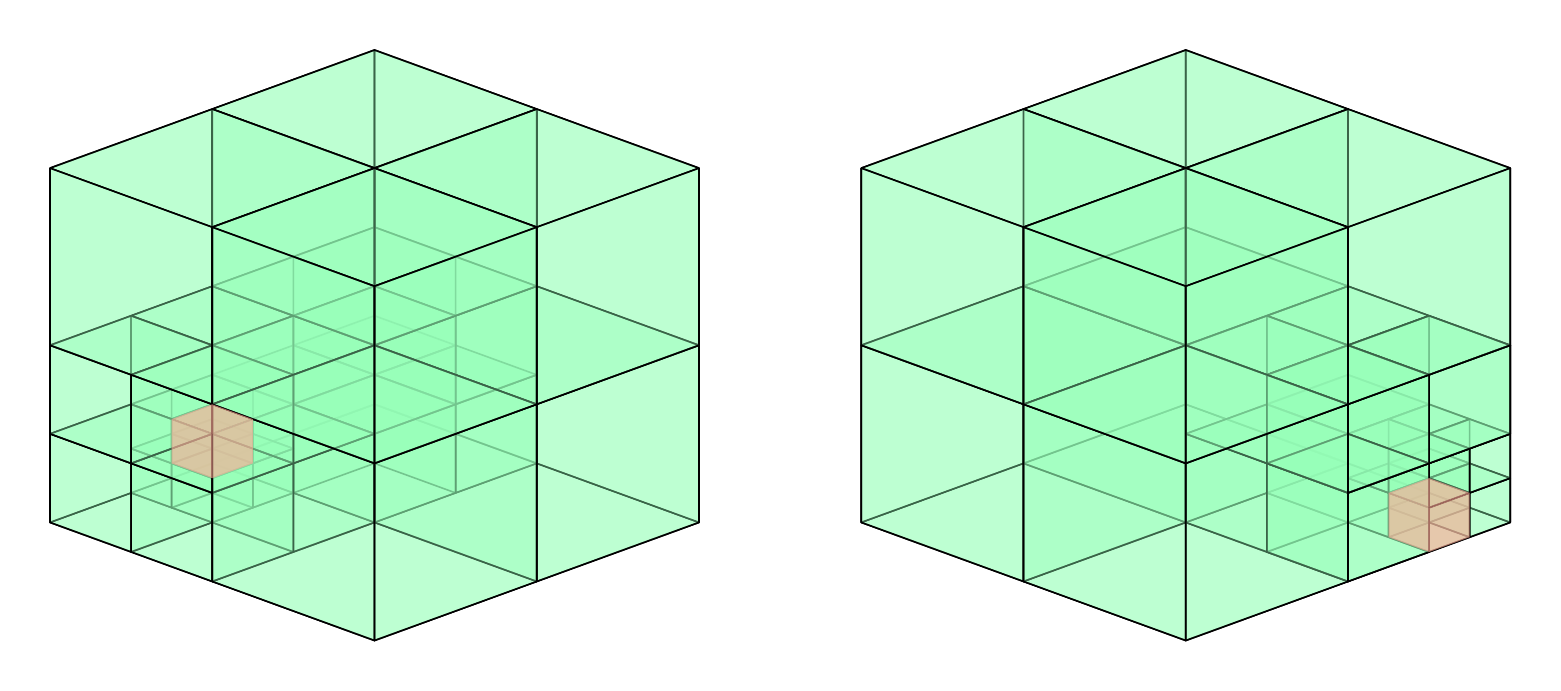}
      \caption{Examples of $h$-adaptive singularity meshes in 3D. Elements containing singularities marked in red.}
    \end{subfigure}

  \label{ch4:examplepoint}
\end{figure}

We define as a \textit{point singularity mesh} a mesh refined hierarchically towards a single point. For example, let us consider a point singularity in a $d$-dimensional space refined towards a point $Q$ until some arbitrary refinement level $r$. We denote such mesh as $\mathcal{G}_Q^d(r)$. The singularity point can be placed either inside an element or on an element's boundary. In particular, the singularity point can be placed on the whole mesh boundary.  Figure \ref{ch4:examplepoint} illustrates examples of meshes with point singularities.

\begin{figure}[h!]
  \caption[Example point singularity mesh construction]{A $G_q^2(4)$ mesh construction for a point $q$ close to the middle of the mesh.}
  \includegraphics[width=\textwidth]{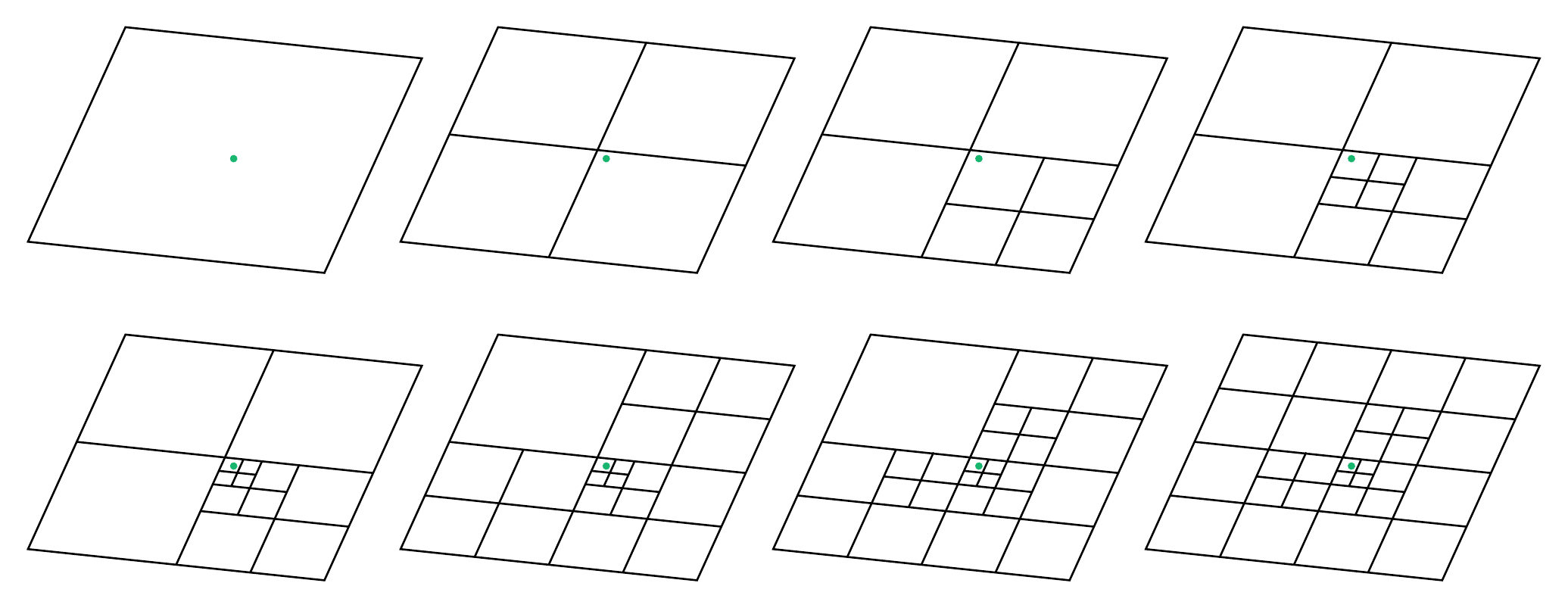}
  \centering
  \label{ch4:exampleconstruct}
\end{figure}

Figure \ref{ch4:exampleconstruct} illustrates a process of building such a mesh using Algorithm \ref{algo:mesh-construction}.

\subsection{Analysis: point singularity placed on the boundary of the mesh}
\begin{figure}[h]
\includegraphics[width=0.7\textwidth]{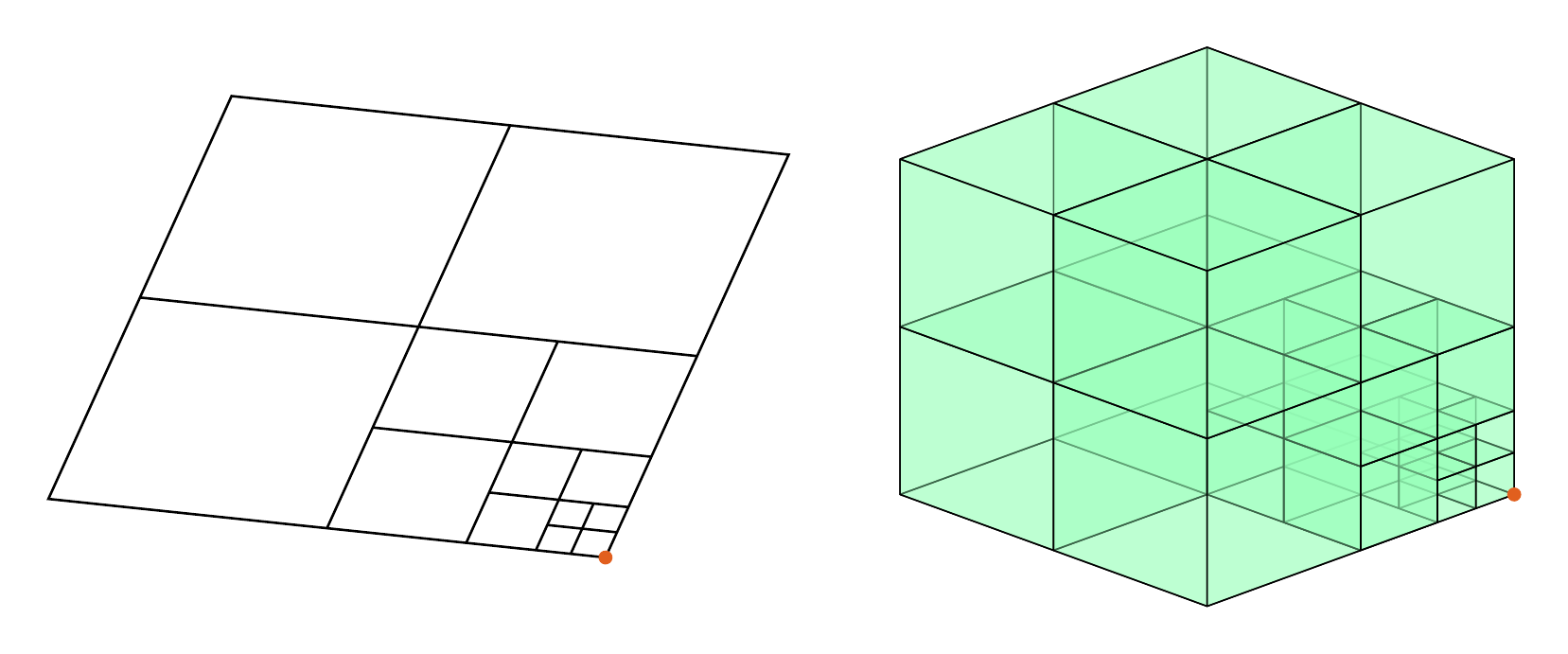}
\centering
\caption[Meshes in 2D and 3D with corner point singularity.]{Meshes in 2D and 3D with corner point singularity marked in red.}
\label{ch4:cornerexample}
\end{figure}

\begin{figure}[h]
  \centering
  \includegraphics[width=0.55\textwidth]{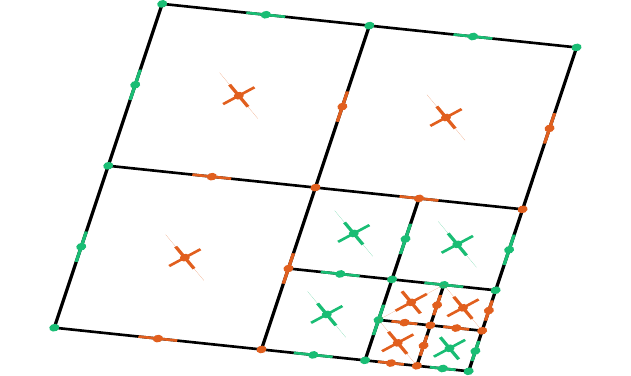}
  \caption[Number of nodes in corner singularity.]{Each layer of refinement adds $2^d - 1$ elements and each element adds $p^d$ \Aniarev{variables}.}
  \label{ch4:cornernodes}
\end{figure}

As an example of a point singularity mesh,  we analyze a relatively simple case with a singularity point $Q$ placed in one of the corners of the mesh. We show a way to generate an ordering that results in the linear complexity of the solver. Figure \ref{ch4:cornerexample} shows examples of $h$-adaptive meshes of this type in two and three dimensions.
For this kind of meshes, it is not required to perform an extra step to ensure that the $1$-irregularity rule is met. Therefore, the number of mesh elements $N_e$ is described by the following formula:
\begin{equation}
N_e(G_Q^d(r)) =
 \begin{cases}
   1 & \text{if } r = 0, \\
   N_e(G_Q^d(r-1)) + 2^d - 1 & \text{if } r > 0,
  \end{cases}
\label{ch4:form:elements1}
\end{equation}
where $r$ is the refinement level of the mesh,  and $d$ is the dimensionality of the mesh. Formula \eqref{ch4:form:elements1} is equivalent to:
\begin{align}
N_e(G_Q^d(r)) &= r(2^d - 1) + 1= \mathcal{O}(2^d r) .
\end{align}
In other words,  the number of elements in a corner point singularity mesh grows linearly with the refinement level $r$.\\
After the first refinement, the mesh has $(2p+1)^d$ \Aniarev{variables}. Further refinement level adds a layer of $(2^d-1)$ elements, and for each new element, $p^d$ \Aniarev{variabes} are created. The following formula describes the number of the basis functions on such a mesh for $r\geq1$:
\begin{equation}
\displaystyle N_v(G_Q^d(r)) = {(2p+1)}^d + (r-1)(2^d - 1) p^d = \mathcal{O}(2^drp^d).
\label{ch4:cornerVarsCount}
\end{equation}
If $p$ is constant, the number of variables is linearly proportional to the refinement level $r$. Figure \ref{ch4:cornernodes} visually explains how the number of basis functions grows.
\subsubsection{Time complexity}
\begin{figure}[!htb]
  \centering
  \def\svgwidth{.9\textwidth}
  \includegraphics[width=0.75\textwidth]{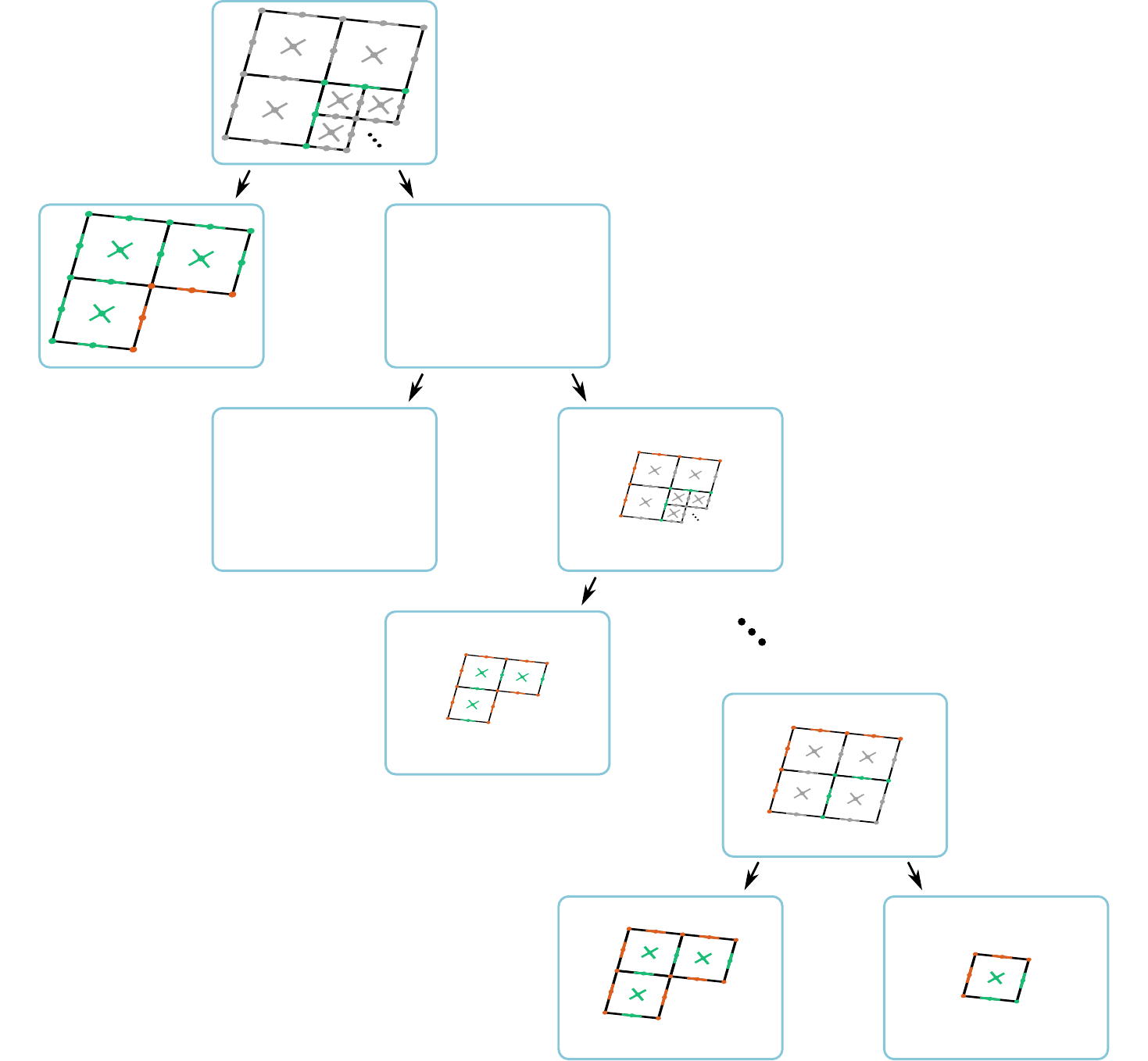}
  \caption[Corner point singularity element partition tree in 2D.]{Corner point singularity element partition tree in 2D. Green nodes are removed at given partition tree node, gray nodes are removed deeper in the tree, and red nodes are removed in the ancestor of the tree node.}
  \label{ch4:cornerEPT2D}
\end{figure}

We use the element partition tree method of ordering generation to calculate the time complexity. The element partition tree can be built by recursively removing the layer of $2^d-1$ least refined (or largest) elements.  Figure \ref{ch4:cornerEPT2D} illustrates the structure of such an element partition tree.
Let us denote each layer by $s_i$, where $i$ is the number of the shell, counting from the most refined one, $s_0$ being the particular case of the single element in the corner. It is easy to see that the whole tree is recursive and that each layer $s_i$, with $i > 0$, contains $2^d - 1$ elements.\\
During the solution process at each layer $s_i$ ($i \geq 1$), we first remove all the rows corresponding to the nodes with support entirely within $s_i$ (denoted by $\mathit{If}(s_i)$). Then, we remove all rows that correspond to nodes on the interface between $s_{i-1}$ and $s_i$ (denoted $\mathit{If}(s_{i-1}, s_i)$). This results in the following complexity formula for a single layer (excluding layers $s_0$, $s_1$ and $s_r$):
\begin{equation}
\begin{array}{r @{\;} l}
  T(s_i) &=  C_r(|\mathit{If}(s_i)|, |\mathit{If}(s_i)| + |\mathit{If}(s_{i-1}, s_i)| + |\mathit{If}(s_{i}, s_i+1)|) \smallskip \\
          &+  C_r(|\mathit{If}(s_{i-1}, s_i)|, |\mathit{If}(s_{i-1}, s_i)| + |\mathit{If}(s_{i}, s_i+1)|) \smallskip \\
        & =  |\mathit{If}(s_i)|\left(|\mathit{If}(s_i)| + 2|\mathit{If}(s_{i-1}, s_i)|\right)^2 + 4|\mathit{If}(s_{i-1}, s_i)|^3.
\end{array}
\end{equation}
To calculate the number of \Aniarev{variables} on the interface between layers,  we realize that this number stays the same from when the given layer is refined. In other words, the number of \Aniarev{variables} on the interface between any layers is equal to the number of \Aniarev{variables} on the interface between layers $s_0$ and $s_1$:
\begin{equation}
|\mathit{If}(s_{i-1}, s_i)| = {(p+1)}^d - p^d.
\end{equation}
Similarly, we can calculate the number of "internal" \Aniarev{variables},  which is equal to the number of \Aniarev{variables} in a mesh of refinement level 1 ($G^d(1)$), excluding the \Aniarev{variables} on the (hyper)faces not touching the corner with singularity and minus the number of \Aniarev{variables} in a mesh of refinement level 0 ($G^d(0)$):
\begin{equation}
|\mathit{If}(s_i)| = (2^d-1)p^d - |\mathit{If}(s_{i-1}, s_i)| = (2^d-1)p^d + p^d - (p+1)^d = (2p)^d - (p+1)^d
\end{equation}
Therefore, the complexity of the elimination of a single layer $s_i$, where $1 < i < r$, goes as follows:
\begin{equation}
\begin{array}{r @{\;} l}
  T(s_i) & = |\mathit{If}(s_i)|(|\mathit{If}(s_i)| + 2|\mathit{If}(s_{i-1}, s_i)|)^2 + 4|\mathit{If}(s_{i-1}, s_i)|^3 \smallskip \\
  & = ((2p)^d - (p+1)^d)((2p)^d - {(p+1)}^d + {(p+1)}^d - p^d)^2 + 4({(p+1)}^d - p^d)^3 \smallskip \\
  & = ((2p)^d - (p+1)^d)((2^d-1)p^d)^2 + 4({(p+1)}^d - p^d)^3 \smallskip \\
  & \leq (2p)^d((2p)^d)^2 + 4(p+1)^{3d} = (2p)^{3d} + 4(p+1)^{3d} = \mathcal{O}{\left((2p)^{3d}\right)}.
  \end{array}
\end{equation}
The computational complexity of the whole mesh is then given as follows:
\begin{align}
 T(G^d(r)) = \mathcal{O}(r(2p)^{3d}) = \mathcal{O}\left((2p)^{2d} \cdot 2^drp^d\right) \mathcal{O}\left(2^{2d}p^{2d}N \right).
\end{align}
With the assumption that the parameter $p$ is constant, and considering that dimension $d$ is constant for a given problem, we deduce that the time complexity of the whole algorithm is linear $\mathcal{O}(N)$.
\subsection{Quasi-optimal \Aniarev{h-adapted mesh towards a point singularity}}\label{ch4:sec:quasioptimal}
To generalize the analysis above,  we present a set of properties of the mesh that guarantee that it is possible to create an ordering that leads to linear solver execution time.
 We call a mesh that fulfills those properties {\em quasi-optimal point singularity mesh}. The properties of {\em quasi-optimal h-adapted mesh towards a point singularity} are the following:
\begin{enumerate}
  \item The mesh can be split into no more than $Kr$ consecutive layers, where $r$ is the refinement level of the mesh, and $K$ is some constant.
  \item Each basis function can be assigned to one of the layers so that a pair of basis functions have non-overlapping supports if they are more than $M$ shells apart, where $M$ is some constant.
  \item Each layer has no more than $L$ basis functions assigned, where $L$ is some constant;
\end{enumerate}
If those properties are met, it is possible to create an ordering in which the nodes are removed layer by layer. Removal of each successive layer requires no more than $C_r(L, L(M+1)) = O(L^3M^2)$ subtractions, and the whole solver requires no more than $\mathcal{O}(KrL^3M^2)$ operations. The constructed mesh has $N = \mathcal{O}(KrL)$ nodes, so the time complexity of the solver is stated as $\mathcal{O}(N(LM)^2)$. As $K$, $L$, and $M$ are constants, the resulting solution complexity is $\mathcal{O}(N)$. That is, it grows linearly with the number of nodes (which, in turn, increases linearly with the refinement level).
\subsection{Quasi-optimality of $h$-adaptive point singularity meshes}\label{ch4:sec:generalize}
A final step of the proof presented in this section is to show that any $h$-adaptive point singularity mesh created with Algorithm \ref{algo:mesh-construction} is quasi-optimal. It is easy to notice that:
\begin{enumerate}
\item The elements can be grouped into $R$ layers by their refinement level, and each basis function can be assigned to the layer of one of its elements.
\item Because of the {\em $1$-irregularity rule}, two elements sharing a vertex cannot differ by more than two refinement levels, so a single basis function cannot span over two elements with a difference of more than two refinement levels.
\item There are no more than $4^d$ elements of each refinement level. At the same time, each element has no more than $(p+1)^d$ basis functions -- both values are constant.
\end{enumerate}
Those observations lead to the conclusion that any point singularity $h$-adaptive mesh is quasi-optimal:
\begin{enumerate}
\item There are exactly $R$ layers.
\item Nodes that are more than two layers apart never overlap.
\item There are no more than $4^d(p+1)^d$ nodes at each layer.
\end{enumerate}
Extending those observations and seeing that a mesh is refined towards more than one is straightforward. Still, a finite number of point singularities is quasi-optimal in the same way (each layer potentially has the number of basis functions multiplied by the number of singularities).
\section{Computational complexity for \Aniarev{h-adapted meshes towards a multi-dimensional singularity }}\label{ApB}
This Appendix analyzes the time complexity of direct solvers being run over adaptive meshes refined around singularities of simple shapes of higher dimensionality. In particular, we analyze meshes with singularities in the form of lines, planes/faces, and hyperplanes/hyperfaces with three or more dimensions, depending on the dimensionality of the space.
In this case, for the simplicity of derivation,  we ignore the polynomial order $p$ factor.
\subsection{Structure of the mesh with singularity}
\label{ch5:sec:def}
\begin{figure}[h]
\begin{subfigure}[b]{0.5\textwidth}
\includegraphics[width=0.8\textwidth]{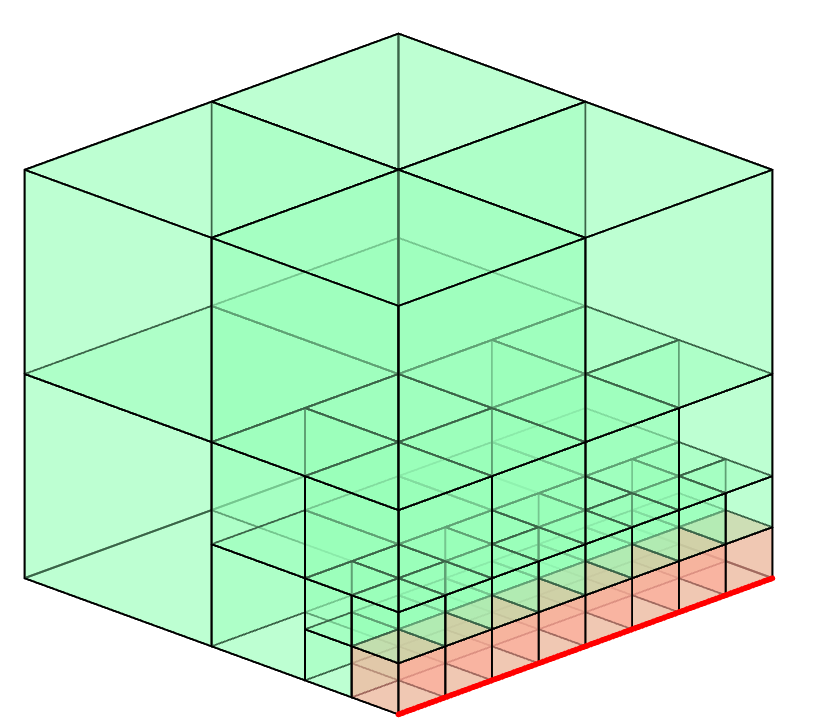}
\centering
\caption{Example of edge singularity mesh construction}
\label{ch5:exampleconstruct3dedge}
\end{subfigure}
\begin{subfigure}[b]{0.5\textwidth}
\includegraphics[width=0.8\textwidth]{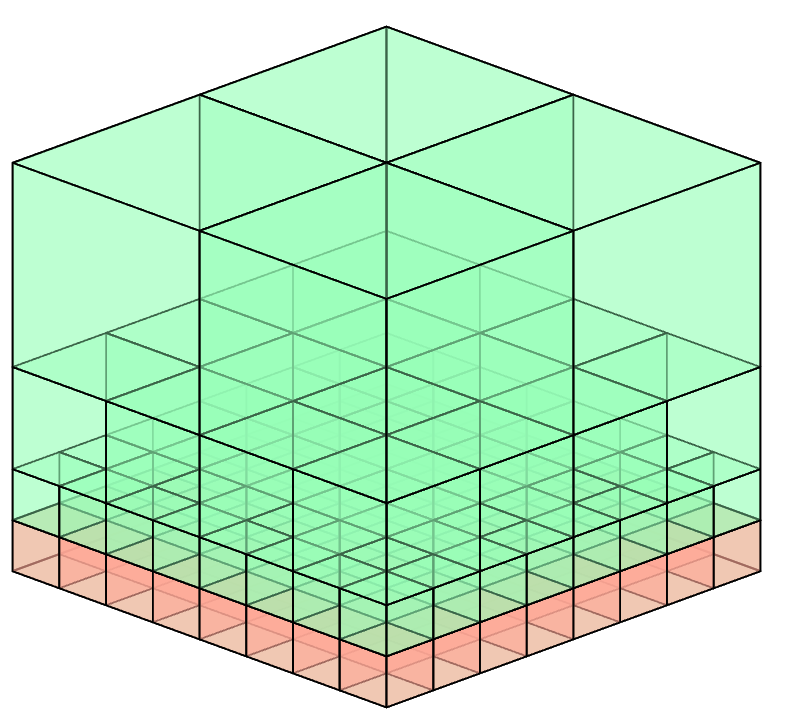}
\centering
\caption{Example of face singularity mesh construction}
\label{ch5:exampleconstruct3dface}
\end{subfigure}
\caption{Examples of mesh singularities in 3D}
\end{figure}

This section analyzes $d$-dimensional meshes refined towards $q$-dimensional singularities. The mesh is considered refined until the refinement level $r$ towards that singularity if all elements overlapping any section containing the singularity have been refined. Such refinement can again be achieved using Algorithm \ref{algo:mesh-construction}. A process of such refinement is shown in Figure \ref{ch5:exampleconstruct3dedge} and
\ref{ch5:exampleconstruct3dface}, respectively
for edge and face singularity in 3-dimensional space.
It is crucial to notice that regular meshes can also be analyzed as mesh refined toward $d$-dimensional singularity in $d$-dimensional space.
\subsection{Analysis: singularity placed on the boundary of the mesh}
For the sake of simplicity,  we start by analyzing the basic case of the singularity placed on the boundary of the full $h$-adapted mesh. Let us denote the dimensionality of the mesh as $d$, the dimensionality of the singularity as $q$, and the singularity itself as $S_q$. The mesh is denoted as $\mathcal{G}_{S_q}^d(r)$, where $r$ is its refinement level.
\subsubsection{Properties of the mesh}
The number of elements in such a singularity mesh can be calculated using the following recursive formula:
\begin{equation}
\begin{array}{r @{\;} l}
N_e(\mathcal G_{S_q}^d(0))& = 1 \smallskip \\
N_e(\mathcal G_{S_q}^d(r)) &= N_v(\mathcal G_{S_q}^d(r-1)) + 2^{q(r-1)}(2^d-1)\smallskip \\
N_e(\mathcal G_{S_q}^d(r)) &=
\begin{cases}
  (2^d-1)r + 1 & \text{if } q = 0 \smallskip\\
  \frac{(2^d-1)2^{rq}-(2^d-2^q)}{2^q-1}  & \text{if } q \geq 1
\end{cases}
\end{array}
\end{equation}
The formula expands to the values given in Table~\ref{ch5:table1}.
\begin{table}[h!]
\centering
\begin{tabular}{ |c|c|c|c|c|c|}
%& \multicolumn{3}{|c|}{Singularity type} \\
 \hline
 & Point ($q=0$) & Edge ($q=1$) & Face ($q=2$) & Hyperface ($q=3$)\\
\hline
1-D 		& $r+1$ 		& $2^r$				& & \\
\hline
2-D  		& $3r+1$ 		& $3(2^r)-2$ 			& $4^r$ 							& \\
\hline
3-D  		& $7r+1$ 		& $7(2^r)-6$ 			& $\tfrac{1}{3}(7(4^r)-4)$  				& $8^r$ \\
\hline
4-D  		& $15r+1$ 		& $15(2^r)-14$ 		& $\tfrac{1}{3}(15(4^r)-12)$  			& $\tfrac{1}{7}(15(8^r) - 8)$ \\
\hline
$d$-D  	& $(2^d-1)r+1$ 	& $(2^d-1)2^r-(2^d-2)$ 	& $\tfrac{1}{3}((2^d-1)4^r - (2^d - 4))$ 		& $\tfrac{1}{7}((2^d-1)8^r - (2^d - 8))$ 	\\
\hline
\end{tabular}
\caption{Number of elements for \Aniarev{d--dimensional h--adapted meshes towards g-dimensional singularities}}
\label{ch5:table1}
\end{table}
For $q \geq 1$, the number of elements is approximated by the following lower and upper bounds:
\begin{equation}
2^{(d-1)}(2^q)^r \leq N_e(\mathcal G_{S_q}^d(r)) \leq 2^d(2^q)^r.
\end{equation}
The number of variables is estimated to be between $p^d$ and $(p+1)^d$ per element. Leading us to the following approximation for $q \geq 1$:
\begin{equation}
2^{(d-1)}p^d(2^q)^r \leq N_v(\mathcal G_{S_q}^d(r)) \leq 2^d(p+1)^d(2^q)^r.
\end{equation}
For set $d$, $q \geq 1$ and $p$, the approximations lead to the following formulas:
\begin{align}
N_e(\mathcal G_{S_q}^d(r)) & = \mathcal{O}((2^q)^r),\\
N_v(\mathcal G_{S_q}^d(r)) & = \mathcal{O}((2^q)^r).
\end{align}
In other words, both the number of elements and variables grow proportionally to $\mathcal{O}(2^{qr})$ and this growth speed (understood in terms of $\mathcal{O}$-notation) depends only on the dimensionality of the singularity $q$, not on the dimensionality of the mesh $d$.
\subsubsection{Time complexity of a solution with singularity built on mesh boundary}
\begin{figure}[ht]
\includegraphics[width=0.75\textwidth]{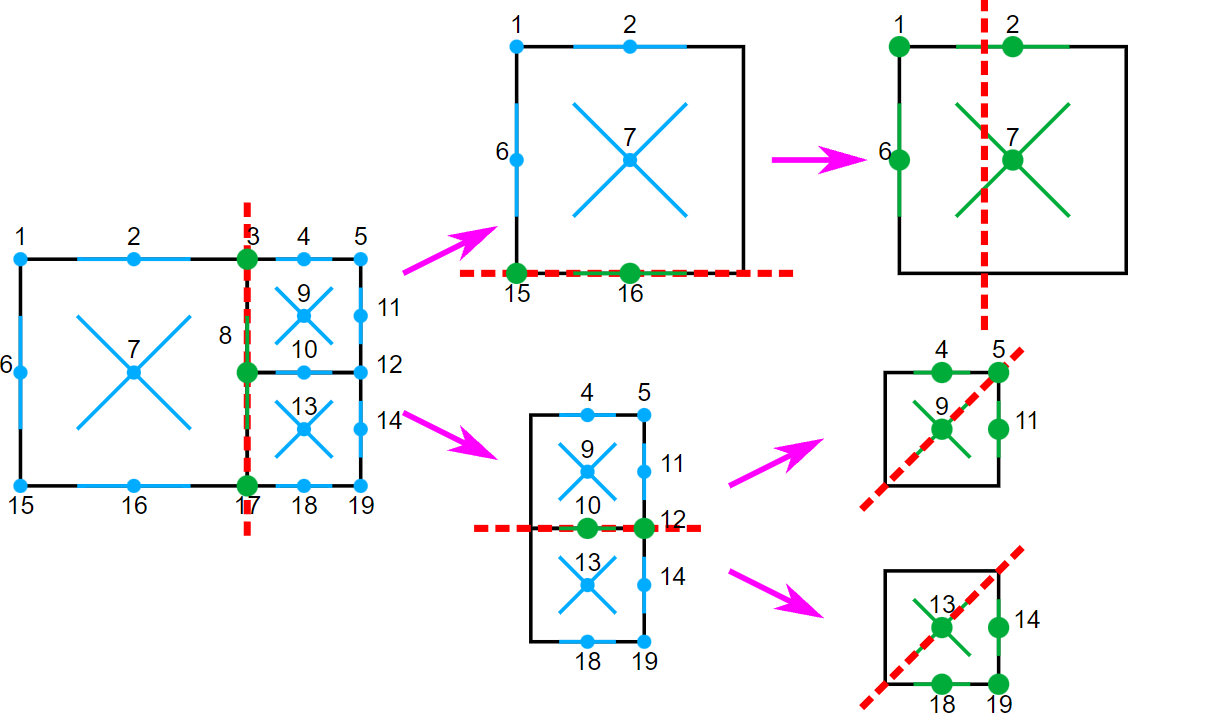}
\centering
\caption{Example element partition tree for one-dimensional boundary singularity in 2-D.}
\label{ch5:partitionExample}
\end{figure}

To analyze the time complexity of the solver, we can again use the element partition tree approach. The element partition tree is built using the following recursive procedure:
\begin{enumerate}
\item Create a root node of the element partition tree and attach all the elements to that node.
\item If there is just one element, finish the procedure, and the root node is the sole node of the returned tree.
\item Create a child node of the root node containing all the least refined elements.
\item Divide the remaining elements by $q$ parallel planes perpendicular to the singularity and parallel to the boundaries of the mesh (let us denote those planes as {\em dividing planes}), crossing the midpoint of the singularity. This refinement creates $2^q$ sub meshes.
\item For each sub mesh generated above, run this procedure recursively and attach the resulting trees as subtrees of the second child node.
\item Finish the procedure and return the tree stemming from the root node.
\end{enumerate}
Figure \ref{ch5:partitionExample} shows an example of such an element partition tree. We denote by $s_i$ the element partition tree nodes from the refinement level $i$.\\
Even though the \Aniarev{element partition} tree nodes on each level have an analogous set of elements, the order of elimination differs slightly for the tree nodes that contain elements on the boundary of the mesh other than the boundary containing the singularity. To simplify the analysis, we modify the order of elimination slightly so that those tree nodes behave similarly to the others: the variables corresponding to the basis functions on the boundary of the mesh are eliminated at the root node $s_r$. This change increases the computation time slightly. However, it has no impact on the time complexity.
To calculate the computational complexity of the solver using the ordering generated from that element partition tree,  we need to know two values for each element partition tree node:
\begin{itemize}
\item The number of variables removed in that element partition tree node $n_r$: For the $s_i$ nodes, this number is proportional to the number of elements of that node that are touching the dividing planes that are used to divide the sub mesh further.
\item The total number of variables with support over the elements in this subtree $n_e$: For the $s_i$ nodes, this number is proportional to the number of elements on the dividing planes of the ancestral tree nodes.
\end{itemize}
It is not difficult to see that the cross-section of the mesh, through the dividing planes, behaves as a mesh in the space of one less dimension built over a singularity of one-dimensionality less than the original one. Thus, for $q > 1$, the following equalities hold true:
\begin{align}
n_r(s_i) & = \mathcal{O}(2^{(q-1)r}q) = \mathcal{O}(2^{(q-1)r}),\\
n_e(s_i) & = \mathcal{O}(2^{(q-1)r}q) = \mathcal{O}(2^{(q-1)r}).
\end{align}
All the remaining nodes have a $\mathcal{O}(1)$ number of elements.
Thanks to those observation, we can calculate the time complexity of running the solver for $q > 1$ using the following equation:
\begin{equation}
\begin{array}{r @{\;} l}
  T(s_0) & = \mathcal{O}(1), \smallskip \\
  T(s_r) & =  \mathcal{O}((2^{3(q-1)})^r), \smallskip\\
  T(\mathcal G_{S_q}^d(r)) & = \mathcal{O}((2^{3(q-1)})^r),
  \end{array}
  \label{ch5:costEquation}
\end{equation}
where the second equality in Equation~\eqref{ch5:costEquation} follows from the fact that:
\begin{equation*}
\begin{array}{r @{\;} l}
T(s_r) &= \displaystyle 2^qT(s_{r-1}) + \mathcal{O}(C_e(n_r, n_r+n_e)) + \mathcal{O}(1) = 2^qT(s_{r-1}) + \mathcal{O}(2^{3(q-1)r})\smallskip \\
&= \displaystyle \mathcal{O}\left(\sum_{h=0}^r 2^{q(r-h)} \cdot 2^{3(q-1)h}\right) = \mathcal{O}\left(\sum_{h=0}^r (2^{q(r-h) + 3(q-1)h})\right) = \mathcal{O}\left(\sum_{h=0}^r 2^{q(r + 2h) - 3h}\right).
\end{array}
\end{equation*}
For $q = 1$,  analogous calculations give us the estimates:
\begin{align}
  n_r(s_i) & = \mathcal{O}(r),\\
  n_e(s_i) & = \mathcal{O}(r).
\end{align}
Moreover,  the time complexity of running the solver follows this equations:
\begin{equation}
\begin{array}{r @{\;} l}
  T(s_0) & = \mathcal{O}(1), \smallskip\\
  T(s_r) & = \mathcal{O}(2^r), \smallskip\\
  T(\mathcal G_{S_1}^d(r)) & = \mathcal{O}(2^r),
\end{array}
\label{ch5:costEquation2}
\end{equation}
where the second equality in Equation~\eqref{ch5:costEquation2} follows from:
\begin{equation*}
\begin{array}{r @{\;} l}
T(s_r)  &=  \displaystyle 2T(s_{r-1}) + \mathcal{O}(C_e(n_r, n_r+n_e)) + \mathcal{O}(1)  = 2T(s_{r-1}) + \mathcal{O}(r^3) \smallskip\\
&= \displaystyle \mathcal{O}\left(\sum_{h=0}^r 2^{r-h} \cdot h^3\right) = \mathcal{O}(2^r \cdot 1^3 + 2^{r-1} \cdot 2^3 + \dots + 2^1 \cdot (r-1)^3 + 1 \cdot r^3).
\end{array}
\end{equation*}
Considering that the number of variables $N_v(\mathcal G_{S_q}^d(r))  = \mathcal{O}((2^q)^r)$,  we estimate the time complexity as a fuction of the number of variables  $q \geq 1$ as:
\begin{equation}
  T(\mathcal G_{S_q}^d(r)) = \mathcal{O}(N_v^{3\frac{q-1}{q}}).
\end{equation}
We sum up the analysis in Table~\ref{ch5:complexitiesTable}.

\begin{table}[h!]
\centering
\begin{tabular}{ |c|c|c|c| }
%& \multicolumn{3}{|c|}{Singularity type} \\
 \hline
 Singularity type                       & Variables         & Operations                    & Operations in $N_v$ \\
\hline
Point 		       & $\mathcal{O}(r)$ 		   & $\mathcal{O}(r)$              & $\mathcal{O}(N_v)$\\
\hline
Edge		       & $\mathcal{O}(2^r)$ 		 & $\mathcal{O}(2^r)$            & $\mathcal{O}(N_v)$\\
\hline
Face  		       & $\mathcal{O}(4^r)$ 		 & $\mathcal{O}(8^r)$	           & $\mathcal{O}({N_v}^{\frac{3}{2}})$ \\
\hline
Hyperface (3-D)  & $\mathcal{O}(8^r)$ 		 & $\mathcal{O}(64^r)$           & $\mathcal{O}({N_v}^2)$ \\
\hline
4-D  & $\mathcal{O}(16^r)$ 		 & $\mathcal{O}(512^r)$           & $\mathcal{O}({N_v}^{2.25})$ \\
\hline
5-D  & $\mathcal{O}(32^r)$ 		 & $\mathcal{O}(4096^r)$           & $\mathcal{O}({N_v}^{2.4})$ \\
\hline
$q$-D & $\mathcal{O}((2^q)^r)$ & $\mathcal{O}((2^{3(q-1)})^r)$ & $\mathcal{O}({N_v}^{3\frac{(q-1)}{q}})$ \\
\hline
\end{tabular}
\caption{\Aniarev{Number of variables and operations for different meshes and singularity dimensions}}
\label{ch5:complexitiesTable}
\end{table}

% \section{Generalized singularity mesh definition}
% A very similar analysis can be used to calculate the complexities of other kinds of meshes with structure based on a singularity.
% ...............................

\subsection{Quasi-optimal \Aniarev{h-adapted meshes towards a $q$ singularity }}
\begin{figure}[!htb]
  \centering
  \includegraphics[width=\textwidth]{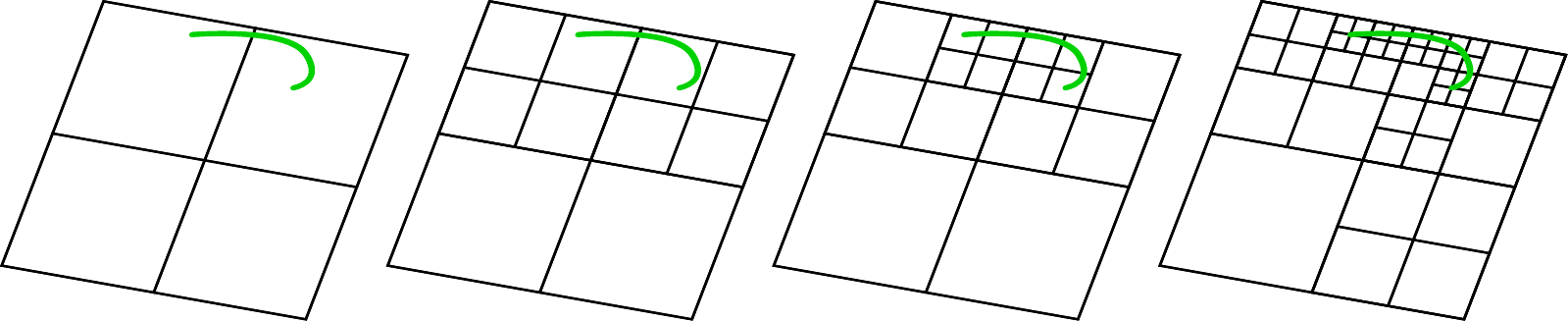}
  \caption[Non-regular edge singularity.]{Two-dimensional mesh refined towards non-regular ``edge'' singularity.}
  \label{ch5:edgenonregular2d}
\end{figure}

\begin{figure}[!htb]
  \centering
  \def\svgwidth{.9\textwidth}
  \includegraphics[width=0.4\textwidth]{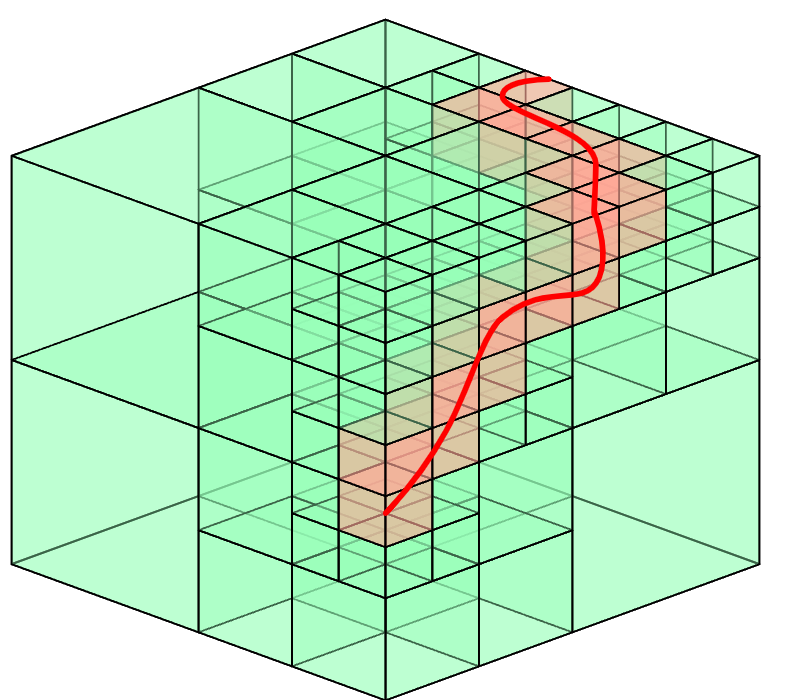}
  \caption[Non-regular edge singularity.]{Three-dimensional mesh refined towards non-regular ``edge'' singularity.}
  \label{ch5:edgenonregular}
\end{figure}

To generalize the analysis from the previous section, we can observe that a broader class of meshes with singularities follow the same time complexity -- we define those meshes as {\em quasi-optimal $q$-dimensional singularity meshes}.
Examples of such meshes are presented in Figures \ref{ch5:edgenonregular2d} and \ref{ch5:edgenonregular}.

We say that a mesh is a quasi-optimal $q$-dimensional singularity mesh if it has the following properties:
\begin{enumerate}
\item The basis functions of the mesh can be assigned to tree nodes of a full $K$-nary tree ($K \geq 2$) of height not larger than some $R = \lceil{\log_K{N}}\rceil+R'$, where $R'$ is some constant.
\item If a pair of basis functions have overlapping supports,  they are assigned to the same tree node,  or one of them is assigned to an ancestor of the tree node of the other one.
\item Each tree node has not more than $Q K^{h\frac{q-1}{q}}$ (if $q \geq 2$) or more than $Q  h$ (if $q \leq 1$) basis functions assigned, where $h$ is the height of the subtree that given node is the root of and $Q$ is some arbitrary positive constant. At the same time, the number of overlaps between basis functions belonging to that tree with basis functions of ancestor tree nodes is limited by the same number.
\end{enumerate}
If such a tree is created,  we can use it to define an ordering that would follow the post-order traversal of that tree. If so, for a $q$-dimensional singularity with \Aniarev{$q \leq 1$}, the complexity of removing all nodes belonging to a tree node $S_h$, with height $h$,  is bounded by the following equation:
\begin{equation}
  T(S_h) \leq C_r(Qh, 2Qh) = 4  Q^3h^3.
\end{equation}
The total complexity of the execution of the solver is in turn no more than:
\begin{equation}
\begin{array}{l @{\;} l}
T(\mathcal G(r)) & \displaystyle \leq \sum_{h=1}^R K^{R-h} \cdot 4(Q^3h^3) = 4Q^3 K^{R}\sum_{i=1}^RK^{-h}h^3\smallskip \\
   & = \displaystyle Q^3 K^{R}\left(\frac{1^3}{K^1} + \frac{2^3}{K^2} + \dots + \frac{(R-1)^3}{K^(R-1)}+ \frac{R^3}{K^R}\right) < Q^3 K^R \cdot 26.
    \end{array}
\end{equation}
Thus,
\begin{equation}
T(\mathcal G(r)) = \mathcal{O}(K^R) = \mathcal{O}(N).
\end{equation}
In the case of $q$-dimensional singularity with $q \geq 2$,  the complexity of removing all nodes belonging to the tree node $S_h$, with height $h$, is bounded by the following equation:
\begin{equation}
  T(S_h) \leq C_r(Q\cdot K^{h\frac{q-1}{q}}, 2 Q\cdot K^{h\frac{q-1}{q}}) = 4Q^3 \cdot K^{3h\frac{q-1}{q}}
\end{equation}
Therefore,  the total complexity of running the whole solver is no more than:
\begin{equation}
\begin{array}{l @{\;} l}
  T(\mathcal G(r)) & \displaystyle = \sum_{h=1}^R K^{R-h}  \cdot 4Q^3(K^{3h\frac{q-1}{q}}) = 4Q^3 \sum_{h=1}^R K^{R - h + 3h\frac{q-1}{q}} \smallskip\\
  &= \displaystyle 4Q^3 \sum_{h=1}^R K^{R}K^{h\frac{3(q-1)-q}{q}}  = 4Q^3 \sum_{h=1}^R K^{R}(K^{\frac{2q-3}{q}})^h\smallskip \\
   & = \displaystyle 4Q^3 K^{R} \sum_{h=1}^R (K^{\frac{2q-3}{q}})^h  < 4Q^3 K^{R}(K^{\frac{2q-3}{q}})^{R+1}.
\end{array}
\end{equation}
Thus,
\begin{equation}
\begin{array}{l @{\;} l}
T(\mathcal G(r)) & = \displaystyle \mathcal{O}(K^{R}(K^\frac{2q-3}{q})^R)
   = \mathcal{O}(K^{R + \frac{2q-3}{q}R})
   = \mathcal{O}(K^{(1 + \frac{2q-3}{q})R})
   = \mathcal{O}(K^{3\frac{q-1}{q}R}) \\
   &= \displaystyle \mathcal{O}(K^{ 3\frac{q-1}{q}\log_{K}{N}})
   = \mathcal{O}((K^{\log_{K}{N}})^{3\frac{q-1}{q}})
   = \mathcal{O}(N^{3\frac{q-1}{q}}).
   \end{array}
\end{equation}
Summarizing,  both cases can be stated as:
\begin{equation}
 \displaystyle T(\mathcal G(r)) = \mathcal{O}\left(N^{\min\left\{3\frac{q-1}{q}, 1\right\}}\right).
\end{equation}
\subsection{Quasi-optimality of $h$-adaptive meshes around singularities}\label{ch5:sec:generalize}
To prove that every $h$-adaptive mesh adapted towards a singularity using the Algorithm \ref{algo:mesh-construction} is quasi-optimal, we propose the following tree generation algorithm:
\begin{enumerate}
  \item Find a dividing plane that crosses the least amount of basis functions' supports out of all planes perpendicular to the singularity that divide the mesh so that no more than half of all basis functions lay solely on either one of the sides of the plane. Create a tree node with all basis functions that the chosen plane crosses the support of and remove them from the mesh.
  \item For each side of the plane, take the basis functions on that side and recursively run the algorithm. The resulting trees become subtrees of the node created in the previous procedure.
  \item Return the tree rooted in the node created in the first step.
\end{enumerate}
%
%
%\begin{figure}[ht]
%  \includegraphics[width=0.75\textwidth]{ph-mid}
%  \centering
%  \caption{A quasi-optimal tree for line singularity.}
%  \label{ch5:quasi-algorithm}
%\end{figure}
%
%
Let us analyze if a tree generated that way proves the quasi-optimality of the mesh:
\begin{enumerate}
\item Every child of any node has at most half of the basis functions of its parent. Because of that, the total height of the tree cannot be larger than $\lceil{\log_2{N_v}}\rceil$ ($K = 2$).
\item There is no overlap between supports of basis functions of either side of the dividing plane.
\item Elements of refinement level $r$ cannot be produced farther than $2\sqrt{d}$ times the side of those elements. Thus the number of elements of a given refinement level crossing the dividing plane is limited.
\begin{enumerate}
\item In the case of $q = 1$, all the dividing planes are parallel to each other. Furthermore, the dividing plane of tree node $g$ nodes deep from the root has a distance between the nearest planes of tree nodes higher in the gree of at most $2^{-g} \cdot L$, where $L$ is the length of the singularity. Last means that no elements larger than $2^{-g} \cdot L$ are eliminated at this node. In other words,  while each dividing plane crosses at most some $Q$ elements of each refinement level ($Q$ is an arbitrary constant), the minimal refinement level of elements removed at given tree node increases by $1$ every step down the tree. Last means that the root node has at most $Q\cdot R$ elements, and other nodes have at most $Q \cdot h$ elements.
\item In case of $q \geq 2$, the limit of elements of refinement level $r$ crossed by a dividing plane is $\mathcal{O}(2^{r(q-1)} \cdot L)$, where $L$ is the length/area/volume of cross section between the singularity and the dividing plane (assuming that the whole mesh has side of length $1$). Because the plane of the smallest cross-section is chosen, $L$ decreases on average by $2^\frac{q-1}{q}$.  The total number of elements crossed are $\mathcal{O}(2^{r(q-1)} \cdot ((\frac{1}{2})^\frac{q-1}{q})^{R-h})$. Considering that $R = rq$ (up to a constant), this is equivalent to $\mathcal{O}(2^{r(q-1)} \cdot ((\frac{1}{2})^\frac{q-1}{q})^{rq-h}) = \mathcal{O}(2^{h\frac{q-1}{q}})$.
\end{enumerate}
\end{enumerate}

\end{document}